\documentclass[reqno]{amsart}

\allowdisplaybreaks
\newtheorem{theorem}{Theorem}[section]
\newtheorem{lemma}[theorem]{Lemma}
\newtheorem{corollary}[theorem]{Corollary}

\theoremstyle{remark}

\numberwithin{equation}{section}

\numberwithin{equation}{section}

\newcommand{ \mr }{ \mathbb{R} }

\def\Xint#1{\mathchoice
    {\XXint\displaystyle\textstyle{#1}}%
     {\XXint\textstyle\scriptstyle{#1}}%
     {\XXint\scriptstyle\scriptscriptstyle{#1}}%
     {\XXint\scriptstyle\scriptscriptstyle{#1}}%
    \!\int}
\def\XXint#1#2#3{{\setbox0=\hbox{$#1{#2#3}{\int}$}
    \vcenter{\hbox{$#2#3$}}\kern-.5\wd0}}

\begin{document}

\title[Very weak solution in Orlicz space]{Gradient estimates of very weak solutions to general quasilinear elliptic equations}

\author[S.-S. Byun, M. Lim]{Sun-Sig Byun\and Minkyu Lim}

\address{Sun-Sig Byun: Seoul National University, Department of Mathematical Sciences and Research Institute of Mathematics, Seoul 151-747, Korea}
\email{byun@snu.ac.kr}

\address{Minkyu Lim: Seoul National University, Department of Mathematical Sciences, Seoul 151-747, Korea}
\email{mk0314@snu.ac.kr}

\thanks{S.  Byun was supported by the National Research Foundation of Korea grant (NRF-2017R1A2B2003877). M. Lim was supported by the National Research Foundation of Korea grant (NRF-2019R1C1C1003844).}



\begin{abstract}
We establish a gradient estimate for a very weak solution to a
quasilinear elliptic equation with  a nonstandard growth condition,
which is a natural generalization of the $p$-Laplace equation. We
investigate the maximum extent for the gradient estimate to hold
without imposing any regularity assumption on the nonlinearity other
than basic structure assumptions. Our results also include a higher
integrability result of the gradient and the existence for the very
weak solutions to such nonlinear problems.
\end{abstract}

\maketitle

\section{Introduction}
In this paper we study the existence and regularity issues regarding
nonlinear elliptic equations with nonstandard growth of the form
\begin{equation}\label{maineq}
\mathrm{div\,}A(x,Du)=\mathrm{div\,} \left(
\frac{\varphi'(|\mathbf{f}|)}{|\mathbf{f}|} \mathbf{f} \right)
\quad  \textrm{in}\ \Omega,
\end{equation}
where $\Omega \in \mr^n$ $(n\geq 2)$ is a bounded open set and
$\varphi \in C^1([0,\infty)) \cap C^2((0,\infty))$ is a given Young
function which is convex and increasing.  A Carath\'eodory map $A =
A(x,\xi): \Omega \times \mr^n \rightarrow \mr^n $ is the
nonlinearity which satisfies the following growth and monotonicity
conditions:
\begin{equation}\label{moncondi}
\begin{cases}
|A(x,\xi)| \leq L\varphi'(|\xi|)\\
\langle A(x,\xi)-A(x,\zeta),\xi-\zeta\rangle \geq \nu
\varphi''(|\xi|+|\zeta|)|\xi-\zeta|^2
\end{cases}
\end{equation}
for a.e $x\in\Omega$, $\xi, \zeta \in\mr^n\backslash \{0\}$ with
some constants $0 < \nu \leq L < + \infty$. We assume that there
exists a constant $s_{\varphi}\geq 1$ such that
\begin{equation}\label{youngcondi}
\varphi(0)=0, \quad \frac{1}{s_{\varphi}}\leq \frac{t
\varphi''(t)}{\varphi'(t)}\leq s_{\varphi} \quad (t >0).
\end{equation}
We are concerned with gradient estimates of a solution to the
problem \eqref{maineq} with only basic structural assumptions
\eqref{moncondi}-\eqref{youngcondi} but without any regularity
assumption on the nonlinearity $x \mapsto A(x,z)$. More precisely, we
want to prove that there exists a positive constant $\delta_{0}$
depending on $n, s_{\varphi}, \nu$ and $L$ such that the following
implication
\begin{equation}\label{implica}
\varphi ( |\mathbf{f}| ) \in L^{q}_{\mathrm{loc}}(\Omega)
\implies  \varphi ( |Du| ) \in L^{q}_{\mathrm{loc}}(\Omega)
\end{equation}
holds for every $q \in [1-\delta_{0}, 1+\delta_{0}]$.  For $q \in [1, 1+\delta_0]$ with
some small $\delta_0>0$, this follows from the classical theory and
 a higher integrability result for a weak solution as in \cite{Ci, DE, Lb}. In this regard, we mainly consider
the case of $ q \in [1-\delta_0, 1)$ with some small $\delta_0>0$ and there is indeed no such an
existence and regularity theory as far as we are concerned in the
literature.

The condition \eqref{youngcondi} includes the classical case of
$\varphi(t)=t^{p}$ with $p>1$. Then our problem \eqref{maineq} is
reduced to the elliptic problems of the $p$-Laplacian type and
implication \eqref{implica} of interest to many researchers becomes
\begin{equation}
\label{implica2} |\mathbf{f}|^p \in L^{q}_{\mathrm{loc}}(\Omega)
\implies   |Du|^p  \in L^{q}_{\mathrm{loc}}(\Omega)
\end{equation}
for some positive number $q$ lying in the range of $\left(
\frac{p-1}{p}, \infty\right)$ depending on further regularity
assumptions on the nonlinearity $A=A(x,z)$. This kind of estimate was first proved by
Calder\'on and Zygmund in \cite{CZ} for the linear case that $p=2$ in
the range of $q \in  (\frac{1}{2}, \infty)$. For the nonlinear case that $p
\neq 2$ and $q \in (1, \infty)$, there have been extensive
regularity results according to given analytic and geometric
settings as in \cite{BW, Ca, Iw, KZ}. On the other hand, the
implication \eqref{implica2} with $q  \in \left( \frac{p-1}{p},
1\right) $ for the $p$-Laplace problem still remains a wide open
problem, partly due to the fact that the duality argument cannot be
applicable to the nonlinear problem, as pointed out by Iwaniec
\cite{Iw}.

There have been notable works with the gradient estimates below the
natural exponent, that is when $q<1$, when the right hand side  of
\eqref{maineq} is a bounded Radon measure, see for instance
\cite{BG1, BG2, Mi, Ph}. Recently, these kinds of estimates have
been extended to the generalized $p$-Laplacian equations in
\cite{Ch2, CGZ, CM}. We refer to \cite{Ch, Mi2} for a further
discussion in this direction. On the other hand, when we have the
divergence data of $\textrm{div} \left( |\mathbf{f}|^{p-2}
\mathbf{f} \right)$ in the right hand side, it was shown in
\cite{IS} that there exists a small positive constant $\delta_{0}$
independent of $\mathbf{f}$ and $u$ such that \eqref{implica2} holds
in the range of $q \in [1-\delta_{0}, 1]$, see also \cite{Ad, IS,
Le}. However, such results have not yet been extended to the
nonlinear equations with nonstandard growth including
\eqref{maineq}. The reason is that a generalization of Hodge
decomposition or Lipschitz truncation method to Orlicz spaces should be properly formalized in advance to deal with very weak
solutions. In this paper we exploit the Lipschitz truncation method
recently developed in \cite{DSSV, DSV} in order to find an
appropriate form of the sub-natural gradient estimates. In this
regard, the present paper provides a window to obtain the
sub-natural estimates for other highly nonlinear problems.

To discuss the regularity of solutions to the equation
\eqref{maineq} with $\mathbf{f}$ having a low degree of
integrability, a few concepts of solutions have been introduced in
the literature. Among them, the very weak solution presented therein
might not have a finite natural energy, and so it requires a test
function to be at least  Lipschitz continuous  in the weak
formulation of the equation \eqref{maineq}. Indeed, various kind of
Lipschitz truncation methods have been developed as in \cite{Ad,
KL, Le}. Here we introduce a suitable Lipschitz truncation
method keeping zero boundary condition for the generalized
$p$-Laplacian equation.

Our paper is organized as follows. In the next section we introduce
the notion of very weak solution to the problem and an
Orlicz-Sobolev space to state our main theorem. In Section
\ref{Sec3} we provide a variety of  analytic tools including a
Lipschitz truncation method which will be commonly used in the later
sections. Existence and higher integrability issues for a very weak
solution to the homogeneous equation will be investigated in Section
\ref{Sec4}. The last section is devoted to dealing with the
comparison estimates in the balls under consideration for proving
Theorem \ref{mainthm}.

\section{Notations and main results}
\label{Sec2} We start this section with the concept of very weak
solution. $u\in W^{1,1}(\Omega)$ is said to be a very weak solution
to \eqref{maineq} if $\varphi'(|Du|) \in  L^{1}(\Omega)$,
$\varphi'(\mathbf{|f|}) \in  L^{1}(\Omega)$ and
\begin{equation}\label{maineq2}
 \int_{\Omega} \big\langle  A(x,Dw), D\varphi \big\rangle \,  dx = \int_{\Omega} \left<  \frac{\varphi'(|\mathbf{f}|)}{|\mathbf{f}|} \mathbf{f} , D\eta  \right>  \,  dx
\end{equation}
holds for every $\eta \in C_{0}^{\infty}(\Omega)$, where we have
denoted by $C_{0}^{\infty}(\Omega)$ to mean that the set of smooth
functions with compact support in $\Omega$.

We next introduce standard notations. $B_\rho(y)$ is the open ball
with center $ y \in \mr^n$ and radius $\rho>0$. For an integrable
function $v$ defined on a bounded measurable subset $E \subset
\mathbb{R}^n$, we briefly denote the integral average over $E$ as
$$
\overline{v}_E:=  \hspace{0.25em} \Xint-_E v(x)\; dx =
\frac{1}{|E|}\int_E v(x)\; dx,
$$
where $|E|$ is the Lebesgue measure of $E$. We now introduce
auxiliary vector fields defined by
\begin{equation}\label{vtdef}
 V(z):=\left[ \frac{\varphi'(|z|)}{|z|} \right]^{\frac{1}{2}}z  \quad \left( z \in\mr^n\backslash \{0\} \right).
\end{equation}
Note that if \eqref{youngcondi} is satisfied, then there holds
\begin{equation}\label{vtcomp}
  |V(z_{1})-V(z_{2})|^2  \approx \varphi''(|z_{1}|+|z_{2}|)|z_{1}-z_{2}|^2 \quad  \left( z_{1} , z_{2} \in\mr^n\backslash \{0\} \right),
\end{equation}
where the implied constant depends on $n$ and $s_{\varphi}$. Then
\eqref{moncondi} yields
\begin{equation}\label{monot}
 \big\langle  A(x,z_{1})-A(x,z_{2}) , z_{1}-z_{2} \big\rangle \geq c |V(z_{1})-V(z_{2})|^2 \quad  \left( z_{1} , z_{2} \in\mr^n\backslash \{0\} \right)
\end{equation}
for some positive constant $c$ depending on $n, \nu$ and
$s_{\varphi}$. Moreover, for any $\varepsilon>0$, we have
\begin{equation}\label{phitvt}
 \varphi(|z_{1}-z_{2}|) \leq  c_{\varepsilon} |V(z_{1})-V(z_{2})|^2 +\varepsilon \varphi(|z_{2}|) \quad  \left( z_{1} , z_{2} \in\mr^n\backslash \{0\} \right) ,
\end{equation}
where the constant $c_{\varepsilon}$ depends on $s_{\varphi}$ and
$\varepsilon$.  For a further discussion regarding \eqref{vtdef}
through \eqref{phitvt}, we refer to  \cite{DE, DSV, Ve}.

As mentioned earlier, this paper aims at obtaining gradient
estimates \eqref{finalesti} below the natural exponent for  very
weak solutions to \eqref{maineq} in an appropriate setting of Orlicz
spaces. Since the associated Young function $\varphi$ satisfies the
so-called $\Delta_{2}$-condition and $\nabla_{2}$-condition when
\eqref{youngcondi} holds (see \cite{RR}), the set $L^{
\varphi}(\Omega)$ consisting of all measurable functions $v$ on
$\Omega$ with
$$
\int_{\Omega} \varphi(|v(x)|)\; dx < \infty
$$
becomes a separable reflexive Banach space, where the norm is given
by
$$
\|v\|_{L^{\varphi}(\Omega)}= \inf  \left\{ \lambda > 0 :
\int_{\Omega} \varphi \left( \frac{|v(x)|}{\lambda}\right) \; dx \leq
1 \right\}.
$$
We remark that with \eqref{youngcondi} and $q$ being sufficiently
close to 1, $L^{{\varphi}^{q}}(\Omega)$ also becomes a separable
reflexive Banach space, where we have defined $\varphi^q(t) :=
[\varphi(t)]^{q}$ for $t\geq 0$.  See Lemma \ref{orlicprop}.  The
Orlicz-Sobolev space $W^{1,\varphi}(\Omega)$ is a function space
consisting of all measurable functions $v \in L^{\varphi}(\Omega)$
whose weak derivatives $Dv$ also belong to
$L^{\varphi}(\Omega;\mr^n)$. The norm of $W^{1,\varphi}(\Omega)$ is
given by
$$
\|v\|_{W^{1,\varphi}(\Omega)}=
\|v\|_{L^{\varphi}(\Omega)}+\|Dv\|_{L^{\varphi}(\Omega)} .
$$
Note that $\overline{C^{\infty}(B_{\rho})}=W^{1,\varphi}(B_{\rho})$
for any ball $B_{\rho} \subset \mr^n$, where the completion is taken
with respect to the $W^{1,\varphi}(B_{\rho})$ norm.
$W_{0}^{1,\varphi}(\Omega)$ is defined as the closure of
$C_0^{\infty}(\Omega)$ in $W^{1,\varphi}(\Omega)$. We refer to
\cite{DT} for a further discussion regarding this Orlicz-Sobolev
space.

We now state our main result in this paper.
\begin{theorem}\label{mainthm}
Assume \eqref{moncondi} and \eqref{youngcondi}. Then there exists a
small positive constant $\delta_{0}=\delta_{0}(n,s_{\varphi}, \nu,
L)$ such that for all $\mathbf{f} \in L^{\varphi^{q}}(\Omega;\mr^n)$
with $q \in [1-\delta_{0}, 1+\delta_{0}]$, any very weak solution $u
\in W^{1,\varphi^{1-\delta_{0}}}(\Omega)$ of \eqref{maineq}
satisfies that
$$u \in W_{ \textrm{loc}}^{1,\varphi^{q}}(\Omega)$$
with the estimate
\begin{equation}\label{finalesti}
 \Xint-_{B_{R/2}}[ \varphi(|Du|)]^q \, dx   \leq c \left[  \left(    \hspace{0.25em} \Xint-_{B_{R}} [\varphi(|Du|)]^{1-\delta_{0}} \, dx  \right)^{\frac{q}{1-\delta_{0}}} +   \hspace{0.25em} \Xint-_{B_{R}} [ \varphi(|\mathbf{f}|)]^q \, dx  \right]
\end{equation}
for some positive constant $c$ depending only on  $ n, s_{\varphi},
\nu, L$, whenever $B_{R} \subset \Omega$ .
\end{theorem}
We would like to point out that the estimate \eqref{finalesti}
becomes the standard Calderón-Zygmund type estimate if  we choose
$q$  in the range of $(1, \infty)$, which holds under a suitable
regularity condition on the nonlinearity $A$ in \eqref{maineq} with
respect to $x$,  as follows from  \cite{ BC, Cho, Ve}. Therefore, the
above theorem implies that we can extend the exponent range to $q
\in [1-\delta_{0}, \infty) $ for the gradient estimate, which
includes the estimates below the natural energy level. On the other
hand, the gradient estimate still holds in the range of $q \in
[1-\delta_{0}, 1+\delta_{0}]$ without any such a condition imposed
on the nonlinearity.

\section{Preliminaries}\label{Sec3}

In this section we provide some analytic tools to prove our main estimate \eqref{finalesti}.  Let $\mathcal{M}$ be the Hardy-Littlewood maximal function defined by
\begin{eqnarray}\label{maxidef}
   \mathcal{M}(f)(y) = \sup_{\rho>0}  \hspace{0.25em} \Xint-_{B_{\rho}(y)}  |f(x)| \; dx \quad ( x \in \mr^n )
\end{eqnarray}
for $f \in L^1_{ \textrm{loc}}(\mr^n)$. We recall that an $A_{s}$-weight ($s>1$) is a non-negative function $w \in L^1_{\textrm{loc}}(\mr^n)$ satisfying
\begin{eqnarray}\label{muckdef}
   [w]_{A_s} := \sup_{B_{\rho} \subset \mr^n} \left( \hspace{0.25em} \Xint-_{B_{\rho}} w(x)\; dx \right) \left( \hspace{0.25em} \Xint-_{B_{\rho}} [w(x)]^{\frac{-1}{s-1}} \; dx \right)^{s-1} < \infty.
\end{eqnarray}
Then we say that $w$ belongs to the Muckenhoupt class $A_{s}$ and the quantity $[w]_{A_s}$ is referred to the $A_{s}$-constant of $w$. It is well known that if $w$ is an $A_{s}$-weight, then there exists a constant $c$ depending only on $s$ and $[w]_{A_s}$ such that
\begin{eqnarray}\label{muckdef2}
\int_{B_{\rho}} [ \mathcal{M}(f)(x) ]^{s} w(x)  \, dx  \leq c \int_{B_{\rho}} |f(x)|^{s} w(x)  \, dx    \quad (B_{\rho} \subset \mr^n)
\end{eqnarray}
for every measurable function $f \in L^1_{ \textrm{loc}}(\mr^n)$ with the right hand side of \eqref{muckdef2} being finite \cite[Chapter 5]{St}.
The following lemma includes standard classical theory regarding maximal functions and $A_{s}$ weights.

\begin{lemma}\label{maximucken}   Let $0<\tau<1$. Then for every nonnegative measurable function $f$ such that $f^{1-\delta} \in L^{1}(\mr^n)$ for some $\delta \in (0,\frac{1-\tau}{2}]$, there holds
\begin{eqnarray}\label{maxiesti}
\int_{B_{\rho}}  [ \mathcal{M} (f^{\tau}) (x)]^{\frac{1-\delta}{\tau}} \, dx  \leq c  \int_{B_{\rho}}  [f(x)]^{1-\delta} \, dx
\end{eqnarray}
for some positive constant $c$ depending on $n$ and $\tau$, whenever $B_{\rho} \subset \mr^n$. Moreover, $[\mathcal{M} (f^{\tau})]^{-\frac{\delta}{\tau}}$ is in the Muckenhoupt class $A_{\frac{1}{\tau}}$ with
\begin{eqnarray}\label{muckenbound}
 [ [ \mathcal{M} (f^{\tau}) (x)]^{-\frac{\delta}{\tau}}]_{A_{\frac{1}{\tau}}} \leq c,
\end{eqnarray}
where the constant $c$ depends on $n$ and $\tau$.
\end{lemma}
\begin{proof}
Since $\frac{1-\delta}{\tau}>\frac{1+\tau}{2\tau}>1$,
\eqref{maxiesti} follows from the boundedness of the maximal
function, where the constant $c$ appearing in \eqref{maxiesti} only
depends on  $n$ and $\tau$, not on $\delta$, see \cite{St} for more details
about maximal function. Now, note that if
$\epsilon \leq \frac{1}{2}$, then for any ball $B_{s}(z) \subset
\mr^{n}$, we have
\begin{eqnarray}\label{muckenproof1}
 \hspace{0.25em} \Xint-_{B_{s}(z)} [\mathcal{M} (f^{\tau})]^{\epsilon} \, dx \leq  c  [\mathcal{M} (f^{\tau}) (z)]^{\epsilon} \quad \textrm{a.e.} \  z\in \mr^{n} ,
\end{eqnarray}
where the constant $c$ depends only on $n$ (see for instance \cite{St} p.214). Then since $\frac{\delta}{1-\tau} \leq \frac{1}{2}$, we observe that for any ball $B_{\rho}(y) \subset \mr^n$,
\begin{eqnarray}\label{muckenproof2}
 \notag && \hspace{0.25em} \Xint-_{B_{\rho}(y)}  [\mathcal{M} (f^{\tau})]^{-\frac{\delta}{\tau}} \, dx  \left( \hspace{0.25em} \Xint-_{B_{\rho}(y)} [\mathcal{M} (f^{\tau})]^{\frac{\delta}{1-\tau}} \, dx   \right)^{\frac{1-\tau}{\tau}} \\
 \notag &&\leq  \hspace{0.25em} \Xint-_{B_{\rho}(y)}  [\mathcal{M} (f^{\tau})]^{-\frac{\delta}{\tau}} \, dx  \left(  \inf_{z \in B_{\rho}(y) } 2^{n}  \hspace{0.25em} \Xint-_{B_{2\rho}(z)} [\mathcal{M} (f^{\tau})]^{\frac{\delta}{1-\tau}} \, dx   \right)^{\frac{1-\tau}{\tau}} \\
&& \leq c \left\{ \inf_{z \in B_{\rho}(y)}   \mathcal{M} (f^{\tau}) (z)  \right\}^{-\frac{\delta}{\tau}}  \left\{ \inf_{z \in B_{\rho}(y)} \  \mathcal{M} (f^{\tau}) (z)  \right\}^{\frac{\delta}{\tau}} \leq c
\end{eqnarray}
for some positive constant $c$ depending only on $n$ and $\tau$. Recalling   \eqref{muckdef}, we conclude  the proof of \eqref{muckenbound}.
\end{proof}

We now  provide several properties of the Young function $\varphi$.
\begin{lemma}\label{orlicprop}
Suppose $\varphi$ satisfies \eqref{youngcondi}. Then the following holds: \newline
(a) For any $t>0$, we have
\begin{equation}\label{orlicprop1}
\frac{1}{s_{\varphi}}+1 \leq \frac{t \varphi'(t)}{\varphi(t)}\leq 1+s_{\varphi}.
\end{equation} \newline
(b) For any $0<\lambda \leq 1 $ and  $1 \leq \Lambda < \infty $, we have
\begin{equation}\label{orlicprop2}
\begin{cases}
\lambda^{1+s_{\varphi}} \varphi(t)  \leq \varphi(\lambda t)\leq  \lambda^{(1/s_{\varphi})+1} \varphi(t)  \quad (t >0) \\
 \Lambda^{(1/s_{\varphi})+1} \varphi(t)  \leq \varphi(\Lambda t)\leq   \Lambda^{1+s_{\varphi}} \varphi(t)   \quad (t >0).
\end{cases}
\end{equation}
(c) There exists a constant $\hat{\delta}=\hat{\delta} (s_{\phi})>0$  such that for  all $q \in [1-\hat{\delta},1+\hat{\delta}]$,
\begin{equation}\label{orlicprop3}
\frac{1}{2s_{\varphi}} \leq \frac{t [\varphi^{q}(t)]''}{[\varphi^{q}(t)]'}\leq 2s_{\varphi} \quad (t >0) ,
\end{equation}
and $\varphi^{q}$ is increasing and convex. Moreover, for every $\delta \in (0, \hat{\delta}]$, there holds
\begin{eqnarray}\label{defint}
\int_{0}^{t} [\varphi (s)]^{-\delta}   \, ds  \leq 2 t [\varphi (t)]^{-\delta}  \quad (t \geq 0)  .
\end{eqnarray}
\end{lemma}

\begin{proof}
We first refer to \cite{Lb} for the finding of (a) and (b). A direct
computation gives
\begin{equation}\label{orlicproof1}
 \frac{t [\varphi^{q}(t)]''}{[\varphi^{q}(t)]'}= \frac{t \varphi''(t)}{\varphi'(t)} + (q-1) \frac{t \varphi'(t)}{\varphi(t)}.
\end{equation}
Therefore, taking $\hat{\delta} =\frac{1}{4s_{\varphi}^2} $, we have
\begin{equation}\label{orlicproof2}
\left| \frac{t [\varphi^{q}(t)]''}{[\varphi^{q}(t)]'}- \frac{t \varphi''(t)}{\varphi'(t)} \right| \leq  \left|  (q-1) \frac{t \varphi'(t)}{\varphi(t)}   \right| \leq  2 \hat{\delta} s_{\varphi} \leq  \frac{1}{2s_{\varphi}}.
\end{equation}
 Combining \eqref{orlicproof2} with \eqref{youngcondi}, we get \eqref{orlicprop3}. Differentiating $\varphi^{q}(t)$, we have
\begin{eqnarray}\label{orlicproof3}
&&\hspace{-13mm} \left[ \varphi^{q}(t) \right]' = q  \varphi'(t)  [ \varphi(t)]^{q-1} \geq 0,  \\
\notag&&\hspace{-13mm} \left[ \varphi^{q}(t) \right]'' = q   [ \varphi(t)]^{q-1}  \left[ \varphi''(t) + (q-1) \frac{[ \varphi'(t)]^2}{ \varphi(t)} \right] \\
\label{orlicproof4}
&&\hspace{-8mm}  \geq q    \varphi^{q-1}(t)  \left[ \varphi''(t) - \hat{\delta}  \frac{[ \varphi'(t)]^2}{ \varphi(t)} \right] \geq q    \varphi^{q-1}(t)  \left[ \varphi''(t) -2s_{\varphi} \hat{\delta}  \frac{ \varphi'(t)}{ t} \right] \geq 0.
\end{eqnarray}
Then we conclude that $\varphi^{q}(t)$ is increasing and convex. To show \eqref{defint}, we first note that (b) implies that $\lim_{t \rightarrow 0}  t [\varphi (t)]^{-\delta} =0$. Therefore, we obtain
\begin{eqnarray}\label{orlicproof4}
\notag && \hspace{-14mm} t [\varphi (t)]^{-\delta} =\int_{0}^{t}  \left\{ s [\varphi (s)]^{-\delta} \right\}' ds \\
\notag && = \int_{0}^{t}  [\varphi (s)]^{-\delta}-\delta \left\{ s  \varphi' (s) [\varphi (s)]^{-1-\delta} \right\} ds \\
 && \hspace{-2mm} \underset{\eqref{orlicprop1}}{\geq} \int_{0}^{t}  [\varphi (s)]^{-\delta}-\frac{1}{2s_{\varphi}}   [\varphi (s)]^{-\delta} ds  \geq \frac{1}{2}\int_{0}^{t}  [\varphi (s)]^{-\delta} ds.
\end{eqnarray}

\end{proof}

Lemma \ref{orlicprop} shows that for any $q \in [1-\hat{\delta},1)$,
the spaces $L^{\varphi^{q}}(\Omega)$ and $W^{1,
\varphi^{q}}(\Omega)$ are well-defined and every function $f \in
L^{\varphi^{q}}(\Omega, \mr^n)$ satisfies $ \varphi' (|f|) \in
L^{1}(\Omega)$. Now we introduce some variation of Young's
inequalities related to the Young function $\varphi$, which will be
frequently used later in Section \ref{Sec4}.

\begin{lemma}\label{youngineq}
Suppose $\varphi$ satisfies \eqref{youngcondi} and let
$\hat{\delta}$ be given in Lemma \ref{orlicprop}. Then for any
$\varepsilon \in (0,1]$ and any $\delta \in (0, \hat{\delta}]$,
there hold
\begin{equation}\label{youngineq1}
 t [\varphi(s)]^{-\delta} \varphi'(s) \leq  \varepsilon [\varphi(t)]^{1-\delta} + c_{\varepsilon} [\varphi(s)]^{1-\delta}  \quad (t, s \geq 0)
\end{equation}
and
\begin{equation}\label{youngineq2}
t [\varphi(t)]^{-\delta} \varphi'(s)\leq  \varepsilon [\varphi(t)]^{1-\delta} + c_{\varepsilon} [\varphi(s)]^{1-\delta}  \quad (t, s \geq 0) ,
\end{equation}
where the constant $c_{\varepsilon}$ depends only on $s_{\varphi}$ and $\varepsilon$.

\end{lemma}
\pagebreak

 \begin{proof}
First observe that for any $\tilde{s} > 0$,
\begin{eqnarray}\label{youngproof}
 \left[ \frac{[\varphi(\tilde{s})]^{1-\delta}}{\tilde{s}} \right]' = \left[ (1-\delta)  \frac{\varphi'(\tilde{s})}{\tilde{s}} - \frac{\varphi(\tilde{s})}{\tilde{s}^{2}} \right] [ \varphi(\tilde{s})]^{-\delta} \geq 0.
\end{eqnarray}
Then since $\frac{[\varphi(\tilde{s})]^{1-\delta}}{\tilde{s}} $ is
increasing, we have
\begin{eqnarray}\label{youngproof2}
 \tilde{t} \  \frac{[\varphi(\tilde{s})]^{1-\delta}}{\tilde{s}}\leq  \tilde{t}  \ \frac{[\varphi(\tilde{t})]^{1-\delta}}{\tilde{t}} + \tilde{s} \  \frac{[\varphi(\tilde{s})]^{1-\delta}}{\tilde{s}} =  [\varphi(\tilde{t})]^{1-\delta} +[\varphi(\tilde{s})]^{1-\delta}
\end{eqnarray}
for any $\tilde{t}, \tilde{s} >0$. Putting $\tilde{t}=\varepsilon t$ and $\tilde{s}= \varepsilon^{-2s_{\varphi}}  s$ for $\varepsilon \in (0,1]$, we obtain
\begin{eqnarray}\label{youngproof3}
\notag && \hspace{-14mm}   t \ \frac{ [\varphi(s)]^{1-\delta}}{s} = \varepsilon t \  \frac{\varepsilon^{-2s_{\varphi}-1}[\varphi(s)]^{1-\delta}}{\varepsilon^{-2s_{\varphi}} s}  \\
\notag && \leq \varepsilon t \ \frac{[ \varepsilon^{-2s_{\varphi}-2}\varphi(s)]^{1-\delta}}{\varepsilon^{-2s_{\varphi}} s}  \underset{\eqref{orlicprop2}}{\leq} \varepsilon t \  \frac{ [\varphi(\varepsilon^{-2s_{\varphi}}s)]^{1-\delta}}{\varepsilon^{-2s_{\varphi}} s} \\
 \notag && \leq  [\varphi(\varepsilon t)]^{1-\delta} +[\varphi(\varepsilon^{-2s_{\varphi}} s)]^{1-\delta} \\
 && \hspace{-1.5mm}  \underset{\eqref{orlicprop2}}{\leq} \varepsilon  [\varphi( t)]^{1-\delta} +\varepsilon^{-4s_{\varphi}^2}[\varphi( s)]^{1-\delta} .
\end{eqnarray}
Taking \eqref{orlicprop1} into account, we get the desired
inequality \eqref{youngineq1}. Similarly, \eqref{defint} implies
that $t [\varphi(t)]^{-\delta}$ is an increasing function.
Therefore, we have
\begin{eqnarray}\label{youngproof4}
\tilde{t} [\varphi(\tilde{t})]^{-\delta}  \frac{\varphi(\tilde{s})}{\tilde{s}} \leq  [\varphi(\tilde{t})]^{1-\delta} +[\varphi(\tilde{s})]^{1-\delta}
\end{eqnarray}
for any $\tilde{t}, \tilde{s} >0$. Putting $\tilde{t}=\varepsilon  t$ and $\tilde{s}= \varepsilon^{-\frac{s_{\varphi}}{2}} s$ for $\varepsilon \in (0,1]$, we have
\begin{eqnarray}\label{youngproof5}
\notag && \hspace{-14mm}  t [\varphi(t)]^{-\delta} \frac{\varphi(s)}{s} =   \frac{\varepsilon t}{ \varepsilon^{\frac{1}{2}}[\varphi(t)]^{\delta}} \times \frac{\varepsilon^{-\frac{s_{\varphi}+1}{2}}\varphi(s)}{\varepsilon^{-\frac{s_{\varphi}}{2}} s}  \\
\notag && \leq  \frac{\varepsilon t}{ [\varepsilon^{(1/s_{\phi})+1}\varphi(t)]^{\delta}}\times  \frac{\varepsilon^{-\frac{s_{\varphi}+1}{2}}\varphi(s)}{\varepsilon^{-\frac{s_{\varphi}}{2}} s}  \underset{\eqref{orlicprop2}}{\leq}   \frac{\varepsilon t}{ [ \varphi(\varepsilon t)]^{\delta}} \times \frac{\varphi(\varepsilon^{-\frac{s_{\varphi}}{2}} s)}{\varepsilon^{-\frac{s_{\varphi}}{2}} s}  \\
 \notag && \leq  [\varphi(\varepsilon t)]^{1-\delta} +[\varphi(\varepsilon^{-\frac{s_{\varphi}}{2}} s)]^{1-\delta} \\
 && \hspace{-1.5mm}  \underset{\eqref{orlicprop2}}{\leq} \varepsilon  [\varphi( t)]^{1-\delta} +\varepsilon^{-s_{\varphi}^2}[\varphi( s)]^{1-\delta},
\end{eqnarray}
which implies the inequality \eqref{youngineq2}.
\end{proof}

We next introduce an Orlicz-Sobolev-Poincar\'e type inequality. The
proof can be found in \cite[Theorem 7]{DE}.

\begin{lemma}\label{poincare} Suppose $\varphi$ satisfies \eqref{youngcondi} and let $\hat{\delta}$ be the number given in Lemma \ref{orlicprop}. Then there exists a constant $\theta=\theta(n,s_{\varphi}) \in [1-\hat{\delta},1)$  such that for any $w \in W^{1, \varphi^{\theta}}(B_\rho)$, there holds
\begin{eqnarray}\label{poin1}
   \Xint-_{B_{\rho}} \varphi \left(\frac{|w-w_{\rho}|}{\rho}\right)  \, dx  \leq c \hspace{0.25em} \left(   \hspace{0.25em} \Xint-_{B_{\rho}}  [\varphi(|Dw|)]^{\theta}  \, dx  \right)^{\frac{1}{\theta}}
\end{eqnarray}
for some positive constant $c$ depending only on $n$ and $s_{\varphi}$.
\end{lemma}

One of the application of this lemma is the following variant
regarding the Lipschitz truncation.

\begin{lemma}
\label{liptrun}
Suppose $\varphi$ satisfies \eqref{youngcondi}. For $v  \in W_{0}^{1,\varphi^{\theta}}(B_{\rho})$ with $B_{\rho} \subset \mr^n$ and
 for any $\lambda>0$, we write
$$
 \quad E_{\lambda}:=  \{x \in B_{\rho} :  \left\{ \mathcal{M}([\varphi(|Dv|)]^{\theta})(x) \right\}^{\frac{1}{\theta}} >\lambda \}.
$$
 Then there exist a Lipschitz function $v_\lambda \in W_{0}^{1,\infty} (B_{\rho}) $ and a positive constant $c$ depending on $n$ and $s_{\varphi}$ such that
$$v_\lambda(x)=v(x) \quad \mathrm{and}   \quad  Dv_\lambda(x)=Dv(x) \quad    \quad \mathrm{a.e.} \ x \in B_{\rho} \backslash E_{\lambda}, $$
where the estimate
$\varphi(|Dv_{\lambda}|)(x)\leq c\lambda$
holds for a.e. $ x \in B_{\rho}$.
\end{lemma}

For any function $v$ in Orlicz-Sobolev spaces, the above lemma provides an Lipschitz function $v_{\lambda}$ identical to $v$ except on the set $E_{\lambda}$ which vanishes as $\lambda$ goes to $\infty$. This lemma allows us to choose  $v_{\lambda}$ as a test function for various divergence types of equations, while the set of difference can be controlled by the original function $v$ and $\lambda$. This kind of approximation was first introduced in \cite{AF}.

The case that $\varphi(t) = t^{p}$ of Lemma \ref{liptrun} was proved
in \cite{DKS}. The proof is based on a Whitney type decomposition.
Since the proof in the literature uses Jensen's inequality and
Poincar\'e type inequality, the lemma still holds for our general
$\varphi$ under the condition \eqref{youngcondi}. See also \cite{DSSV}. For  the practical
applications of the Lipschitz truncation method to very weak
solutions, we refer to \cite{Ad, Le}.

We conclude this section with the following technical lemma which
will be used later for deriving the required gradient estimate.

\begin{lemma}\label{algeb}\cite[Lemma 6.1]{Gi}
Let  $\phi : [\frac{R}{2}, R] \rightarrow \mr^n$ be a bounded non-negative function. Suppose that for every choice of $r_{1}$ and $r_{2}$ such that $\frac{R}{2} \leq r_{1} < r_{2} \leq R$, we have
\begin{equation}\label{technical1}
\phi(r_{1}) \leq d \phi(r_{2}) + \frac{A}{(r_{2}-r_{1})^{\beta}}+ B
\end{equation}
for some positive numbers $A, B, \beta >0$ and $d \in (0, 1)$. Then,
\begin{equation}\label{technical2}
\phi(r_{1}) \leq  c \left( \frac{A}{R^{\beta}}+ B \right)
\end{equation}
for some positive constant $c$ depending on $d$ and $\beta$.
Especially, the constant $c$ continuously depends on
$\beta$.
\end{lemma}

\section{Higher integrability and Solvability}
\label{Sec4} We investigate several properties of a very weak
solution to the equation \eqref{maineq}.
 For the sake of convenience, we use the letter $c$ to mean a generic constants depending
on $n,s_{\varphi}, \nu$ and $L$, where the exact values might be different from line to
line. Similarly, $c(\varepsilon)$  represent a generic constant depending on $n,s_{\varphi}, \nu, L$ and $\varepsilon$. We first introduce a higher integrability result of very weak solution to the following homogeneous equation.
\begin{equation}\label{homeq}
\mathrm{div\,}A(x,Dw)=0 \quad  \textrm{in}\ \Omega.
\end{equation}
\pagebreak
\begin{lemma}\label{higher}
Assume \eqref{moncondi} and \eqref{youngcondi}. Then there exists a
small constant \linebreak $\sigma=\sigma(n, s_{\varphi}, \nu, L)$ such that
every very weak solution $w\in W^{1,\varphi^{1-\sigma}}  ( \Omega)$
to the problem \eqref{homeq} belongs to $w\in
W^{1,\varphi^{1+\sigma}}_{\textrm{loc}}  ( \Omega)$. Moreover, we
have the following estimate
\begin{equation}\label{higheresti}
 \left(  \hspace{0.25em} \Xint-_{B_{\rho}} [\varphi(|Dw|)]^{1+\sigma} \, dx  \right)^{\frac{1}{1+\sigma}}  \leq c   \left( \hspace{0.25em}  \Xint-_{B_{2\rho}}[ \varphi(|Dw|)]^{1-\sigma} \, dx \right)^{\frac{1}{1-\sigma}}  
\end{equation}
for some positive constant $c>0$ depending on  $ n,s_{\varphi}, \nu$ and $L$, whenever $B_{2\rho} \subset \Omega$.
\end{lemma}

\begin{proof}
Suppose $ w \in W^{1,\varphi^{1-\delta}}(\Omega)$ is a solution of \eqref{homeq} for some $\delta \in (0,\frac{1-\theta}{2} ]$,
with $\theta$ given in Lemma \ref{poincare}. Fix any ball $B_{2r}(y) \subset B_{2\rho} \subset \Omega$ and choose a cut-off
function $\eta \in C^{\infty}_{0}(B_{2r}(y))$ such that $|D\eta|
\leq \frac{2}{r},$ $0\leq \eta \leq 1$ and $\eta \equiv 1$ on
$\overline{B_{r}(y)}$. Write
$$v :=(w-\overline{w}_{B_{2r}(y) }) \eta  \in W_{0}^{1,\varphi^{1-\delta}}(B_{2r}(y)).$$
Then we apply Lemma \ref{liptrun} to $v$ to find that for any
$\lambda>0$ one has a Lipschitz function $v_\lambda \in
W_{0}^{1,\infty} (B_{2r}(y) )$ such that $v_\lambda=v$ and
$Dv_\lambda=Dv$ a.e. on $ B_{2r}(y) \backslash E_{\lambda}$ and that
$\varphi(|Dv_{\lambda}|)\leq c\lambda$ for a.e. on $ B_{2r}(y)$ for
some positive constant depending on $n$ and $s_{\varphi}$, where
$$
 \quad E_{\lambda}:=  \{x \in B_{2r}(y) :   \left\{ \mathcal{M}([\varphi(|Dv|)]^{\theta})(x) \right\}^{\frac{1}{\theta}} >\lambda \}.
$$
 Since $C^{\infty}_{0}(B_{2\rho})$ is weak-$*$ dense in $(L^{1}(B_{2\rho}))^{*}$, we can take $v_{\lambda}$  as a test function to the equation \eqref{homeq}. Recalling \eqref{moncondi}, we have
\begin{eqnarray}\label{i1esti1}
\notag &&  \hspace{-20mm}  \int_{B_{2r}(y) \backslash E_{\lambda}} \big\langle  A(x,Dw) , Dv_{\lambda}\big\rangle   \, dx = - \int_{ E_{\lambda}} \big\langle  A(x,Dw) , Dv_{\lambda}\big\rangle   \, dx \\
&&\hspace{25mm} \leq c \varphi^{-1}(\lambda)\int_{ E_{\lambda}}  \varphi ' (|Dw|)   \, dx .
\end{eqnarray}
Multiplying both sides of \eqref{i1esti1} by $\lambda^{-(1+\delta)}$  and integrating from $0$ to $\infty$ with respect to $\lambda$, we obtain
\begin{eqnarray}\label{i1esti2}
\notag &&  \hspace{-5mm} I_{1} : = \int_{0}^{\infty} \lambda^{-(1+\delta)} \left[ \int_{B_{2r}(y) \backslash E_{\lambda}} \big\langle  A(x,Dw) , Dv_{\lambda}\big\rangle   \, dx \,   d\lambda \right] \\
\notag && \leq c  \int_{0}^{\infty} \lambda^{-(1+\delta)} \varphi^{-1} (\lambda)\int_{ E_{\lambda}}  \varphi ' (|Dw|)   \, dx    \,   d\lambda   \\
\notag && \leq c \int_{B_{2r}(y)}  \left[ \int_{0}^{g(x)} \lambda^{-(1+\delta)} \varphi^{-1} (\lambda)  \varphi ' (|Dw|)    \,   d\lambda \right]   \, dx  \\
\notag && = c \int_{ B_{2r}(y)}  \left[ \int_{0}^{ \varphi^{-1} (g(x))} [\varphi (s)] ^{-(1+\delta)} s \varphi ' (s)  \varphi ' (|Dw|)    \,   ds  \right]   \, dx  \\
\notag && \leq c \int_{ B_{2r}(y)} \left[  \int_{0}^{ \varphi^{-1} (g(x))} [ \varphi(s)]^{-\delta}   \varphi ' (|Dw|)    \,   ds  \right] \, dx  \\
&& \hspace{-1.5mm}  \underset{\eqref{defint}}{\leq}  c \int_{ B_{2r}(y)}  \varphi^{-1} (g(x)) [ g(x)]^{-\delta} \varphi ' (|Dw|)       \, dx  ,
\end{eqnarray}
where we have defined
\begin{equation}\label{gdef}
g(x):= \left\{ \mathcal{M}([\varphi(|Dv|)]^{\theta})(x) \right\}^{\frac{1}{\theta}}     \quad (x \in B_{2r}(y)).
\end{equation}
 Applying  Lemma \ref{youngineq} to \eqref{i1esti2} with $t=\varphi^{-1} (g(x)) $ and $ s = |Dw| $, we get
\begin{eqnarray}\label{i1esti3}
\notag && \hspace{-4mm} I_{1} \leq c \left( \int_{ B_{2r}(y)}  [g(x)]^{1-\delta}   \, dx  + \int_{ B_{2r}(y)} [ \varphi(|Dw|)]^{1-\delta}   \, dx \right)\\
\notag && \leq c \left( \int_{ B_{2r}(y)}  [\varphi(|Dv|)]^{1-\delta}  \, dx  + \int_{ B_{2r}(y)}  [\varphi(|Dw|)]^{1-\delta}   \, dx \right)\\
\notag && \leq c \left(  \int_{ B_{2r}(y)}  \left[\varphi \left(\frac{|w-\overline{w}_{B_{2r}(y)} |}{r} \right)\right]^{1-\delta}  \, dx  + \int_{ B_{2r}(y)} [ \varphi(|Dw|)]^{1-\delta}    \, dx \right) \\&& \leq c  \int_{ B_{2r}(y)}  [\varphi(|Dw|)]^{1-\delta}  \, dx .
\end{eqnarray}

Next, we estimate $I_{1}$ from below. We split $ B_{2r}(y) \backslash  B_{r}(y)$ as
\begin{equation}\label{split1}
\notag   D_{1} =\{x \in B_{2r}(y) \backslash  B_{r}(y) :   \mathcal{M}([\varphi(|Dv|)]^{\theta})(x) \leq \delta^{\theta} \mathcal{M}([\varphi(|Dw|)]^{\theta}\chi_{B_{2r}(y)})(x)  \}
\end{equation}
and
\begin{equation}\label{split2}
\notag  D_{2} = \{x \in B_{2r}(y) \backslash  B_{r}(y) :   \mathcal{M}([\varphi(|Dv|)]^{\theta})(x) > \delta^{\theta} \mathcal{M}([\varphi(|Dw|)]^{\theta}\chi_{B_{2r}(y)})(x)  \} .
\end{equation}
Since $Dv_{\lambda}=Dv$ a.e. on $B_{2r}(y) \backslash E_{\lambda}$
and $v=w-\overline{w}_{B_{2r}(y)}$ a.e. on $B_{r}(y)$, we have
\begin{eqnarray}\label{i1esti4}
\notag &&   \hspace{-4mm} I_{1}  = \int_{ B_{2r}(y)}  \int_{g(x)}^{\infty} \lambda^{-(1+\delta)} \big\langle  A(x,Dw) , Dv\big\rangle     \,   d\lambda   \, dx\\
\notag && = \frac{1}{\delta} \int_{ B_{2r}(y)}   [{g(x)}]^{-\delta} \big\langle  A(x,Dw) , Dv\big\rangle     \,    dx \\
\notag && \geq \frac{\nu}{\delta} \int_{ B_{r}(y)}   g^{-\delta} \varphi'' (|Dw|) |Dw|^{2}   \,    dx + \frac{1}{\delta} \int_{ D_{1}}   g^{-\delta} \big\langle  A(x,Dw) , Dv\big\rangle      \, dx \\
\notag && \hspace{10mm} + \frac{1}{\delta} \int_{  D_{2}}   g^{-\delta} \big\langle  A(x,Dw) , wD\eta \big\rangle       \, dx \\
\notag && \geq \frac{1}{\delta} \left( \nu \int_{ B_{r}(y)}   g^{-\delta} \varphi'' (|Dw|) |Dw|^{2}  \,    dx - L \int_{ D_{1}}   g^{-\delta}  \varphi' (|Dw|) |Dv|     \, dx  \right. \\
\notag && \hspace{10mm} \left. - 2L \int_{  D_{2}}   g^{-\delta} \varphi' (|Dw|) \frac{|w-\overline{w}_{B_{2r}(y)}|}{r}     \, dx  \right)   \\
 && = :  \frac{1}{\delta}(I_{2} - I_{3} - I_{4}) .
\end{eqnarray}
We first estimate $I_{2}$ from below. According to Lemma \ref{maximucken}, ${g}^{-\delta} \in  A_{1/\theta}$. Then from the boundedness of maximal function(Lemma \ref{maximucken}), we conclude that
\begin{eqnarray}\label{i2esti}
\notag && \hspace{-4mm} I_{2}= \nu \int_{ B_{r}(y)}   g^{-\delta} \varphi'' (|Dw|) |Dw|^{2}  \,    dx    \geq c \int_{ B_{r}(y)}   g^{-\delta} \varphi (|Dw|) \,    dx   \\
&& \geq c \int_{ B_{r}(y)}  g^{-\delta} \left\{ \mathcal{M}([\varphi(|Dw|)]^{\theta}\chi_{B_{r}(y)}) \right\} ^{\frac{1}{\theta}}      \, dx .
\end{eqnarray}
\pagebreak
We now find a pointwise upper bound of $g$. If $x \in B_{r/2}(y)$, then
\begin{eqnarray}\label{gbound}
\notag  &&  \hspace{-8mm} g(x) = \left\{ \mathcal{M}([\varphi(|Dv|)]^{\theta})(x)  \right\}^{\frac{1}{\theta}} \\
\notag && \leq  \sup_{0 < s \leq r - |x-y| } \left(  \hspace{0.25em} \Xint-_{ B_{s}(x)  } [\varphi(|Dv|)]^{\theta}  \, dz \right)^{\frac{1}{\theta}} + \sup_{s > r - |x-y| } \left(  \hspace{0.25em} \Xint-_{ B_{s}(x)  }[ \varphi(|Dv|)]^{\theta}  \, dz \right)^{\frac{1}{\theta}} \\
\notag && \leq   \left\{ \mathcal{M}([\varphi(|Dw|)]^{\theta}\chi_{B_{r}(y)})(x)\right\}^{\frac{1}{\theta}} + c \left(  \hspace{0.25em} \Xint-_{ B_{2r}(y) }[ \varphi(|Dv|)]^{\theta}  \, dz \right)^{\frac{1}{\theta}}\\
&& \leq  \left\{ \mathcal{M}([\varphi(|Dw|)]^{\theta}\chi_{B_{r}(y)})(x)\right\}^{\frac{1}{\theta}}  + c \left(  \hspace{0.25em} \Xint-_{ B_{2r}(y) } [\varphi(|Dw|)]^{\theta}  \, dz \right)^{\frac{1}{\theta}},
\end{eqnarray}
where we have used Lemma \ref{poincare} for the last inequality.  Comparing the last two terms in \eqref{gbound}, we have

\begin{align}\label{algebraic}
\notag  &   \left\{ \mathcal{M}([\varphi(|Dw|)]^{\theta}\chi_{B_{r}(y)})(x) \right\}^{\frac{1-\delta}{\theta}}  \\
& \leq 2  [g(x)]^{-\delta} \left\{ \mathcal{M}([\varphi(|Dw|)]^{\theta}\chi_{B_{r}(y)})(x)\right\}^{\frac{1}{\theta}}  + \left[ c \left(  \hspace{0.25em} \Xint-_{ B_{2r}(y) } [\varphi(|Dw|)]^{\theta}  \, dz \right)^{\frac{1}{\theta}} \right]^{1-\delta}
\end{align}
for  $x \in B_{r/2}(y)$. This estimate leads to
\begin{eqnarray}\label{i2esti2}
\notag && \hspace{-4mm} I_{2} \geq  c \int_{ B_{r/2}(y)}  [g(x)]^{-\delta} \left\{ \mathcal{M}([\varphi(|Dw|)]^{\theta}\chi_{B_{r}(y)})\right\}^{\frac{1}{\theta}}      \, dx  \\
\notag && \geq c \int_{ B_{r/2}(y)} \left\{  \mathcal{M}([\varphi(|Dw|)]^{\theta}\chi_{B_{r}(y)})\right\}^{\frac{1-\delta}{\theta}} \, dx - c r^n \left(  \hspace{0.25em} \Xint-_{ B_{2r}(y) } [\varphi(|Dw|)]^{\theta}  \, dz \right)^{\frac{1-\delta}{\theta}}\\
&& \geq c \int_{ B_{r/2}(y)}  [\varphi(|Dw|)]^{1-\delta}   \, dx - c r^n \left(  \hspace{0.25em} \Xint-_{ B_{2r}(y) } [\varphi(|Dw|)]^{\theta}  \, dx \right)^{\frac{1-\delta}{\theta}}  .
\end{eqnarray}
Using the definition of $D_{1}$ and Lemma \ref{maximucken}, $I_{3}$ can be estimated as follows:
\begin{align}\label{i3esti}
\notag   I_{3} &= L \int_{ D_{1}}   [g(x)]^{-\delta}  \varphi' (|Dw|) |Dv|     \, dx   \\
\notag
&\leq L \int_{ B_{2r}(y)}  \left\{ \mathcal{M}([\varphi(|Dv|)]^{\theta}) \right\}^{-\frac{\delta}{\theta}}  \varphi^{-1}( \left\{ \mathcal{M}([\varphi(|Dv|)]^{\theta})\right\}^{\frac{1}{\theta}}) \varphi' (|Dw|) \, dx \\
\notag  &  \hspace{-1.5mm}  \underset{\eqref{youngineq2}}{\leq} \varepsilon \int_{ B_{2r}(y)}   [\varphi(|Dw|)]^{1-\delta}  \, dx +c(\varepsilon) \int_{ B_{2r}(y)}  \left\{\mathcal{M}([\varphi(|Dv|)]^{\theta})\right\}^{\frac{1-\delta}{\theta}} \, dx  \\
\notag &\leq \varepsilon \int_{ B_{2r}(y)}   [\varphi(|Dw|)]^{1-\delta}  \, dx +c(\varepsilon) \delta^{1-\delta}\int_{ B_{2r}(y)}  \left\{ \mathcal{M}([\varphi(|Dw|)\chi_{B_{2r}(y)}]^{\theta})\right\}^{\frac{1-\delta}{\theta}} \, dx \\
 &\leq ( \varepsilon + c(\varepsilon) \delta^{1-\delta} ) \int_{ B_{2r}(y)}   [\varphi(|Dw|)]^{1-\delta}  \, dx .
\end{align}
Similarly,  we estimate $I_{4}$ using the definition of $D_{2}$. Here we write
\begin{equation}\label{hdef}
 h(x) := \left\{ \mathcal{M}([\varphi(|Dw|)\chi_{B_{2r}(y)}]^{\theta})(x) \right\}^{\frac{1}{\theta}}   \quad (x \in B_{2r}(y))
\end{equation}
 and observe that $\varphi(|Dw(x)|) \leq h(x)$ for $x \in B_{2r}(y)$. Then we find that
\begin{align}\label{i4esti}
\notag  I_{4}& = 2L \int_{ D_{2}}   [g(x)]^{-\delta} \varphi' (|Dw|) \frac{|w-\overline{w}_{B_{2r}(y)}|}{r}   \, dx \\
\notag &  = 2L \int_{ B_{2r}(y)}  \left\{ \mathcal{M}([\varphi(|Dv|)]^{\theta}) \right\}^{-\frac{\delta}{\theta}}  \varphi' (|Dw|) \frac{|w-\overline{w}_{B_{2r}(y)}|}{r} \, dx \\
\notag &  \leq 2 \delta^{-\delta} L  \int_{ B_{2r}(y)}  [h(x)]^{-\delta}  \varphi' ( \varphi^{-1}(h(x)) ) \frac{|w-\overline{w}_{B_{2r}(y)}|}{r} \, dx \\
\notag&  \hspace{-1.5mm}  \underset{\eqref{youngineq1}}{\leq} \varepsilon \int_{ B_{2r}(y)} [h(x)]^{1-\delta}  \, dx  +c(\varepsilon) \int_{ B_{2r}(y)}  \left[\varphi\left(\frac{|w-\overline{w}_{B_{2r}(y)}|}{r}\right)\right]^{1-\delta}  \, dx\\
& \leq  \varepsilon \int_{ B_{2r}(y)}  [\varphi(|Dw|)]^{1-\delta} \, dx + c(\varepsilon) r^n \left(  \hspace{0.25em} \Xint-_{ B_{2r}(y) } [\varphi(|Dw|)]^{\theta}  \, dx \right)^{\frac{1-\delta}{\theta}} .
\end{align}
Combining \eqref{i1esti3}, \eqref{i2esti2}, \eqref{i3esti} and \eqref{i4esti}, we have
\begin{eqnarray}\label{higher1}
\notag && \hspace{-10mm} \int_{ B_{r/2}(y)}  [\varphi(|Dw|)]^{1-\delta}   \, dx \leq  c(\varepsilon) r^n \left(  \hspace{0.25em} \Xint-_{ B_{2r}(y) } [\varphi(|Dw|)]^{\theta}  \, dx \right)^{\frac{1-\delta}{\theta}}  \\
 && \hspace{5mm} + ( c\delta + 2\varepsilon + c(\varepsilon) \delta^{1-\delta}  ) \int_{ B_{2r}(y)}  [\varphi(|Dw|)]^{1-\delta} \, dx .
\end{eqnarray}
Therefore, taking $\varepsilon=\frac{1}{16}$ and $\delta < \tilde{\delta}_{1}:= \min \{\frac{1}{8c},\frac{1}{32c(\varepsilon)^2}, \frac{1-\theta}{2}\}$, we get
\begin{eqnarray}\label{higher2}
\notag && \hspace{-15mm}  \hspace{0.25em} \Xint-_{ B_{r/2}(y)}  [\varphi(|Dw|)]^{1-\delta}   \, dx \leq  c \left(  \hspace{0.25em} \Xint-_{ B_{2r}(y) } [\varphi(|Dw|)]^{\theta}  \, dx \right)^{\frac{1-\delta}{\theta}}  \\
 && \hspace{35mm} + \frac{1}{2} \hspace{0.25em} \Xint-_{ B_{2r}(y)}  [\varphi(|Dw|)]^{1-\delta} \, dx .
\end{eqnarray}
Since the constant $c$ in the above expression does not depends on $\delta$, using the standard Gehring's argument, we get the desired conclusion.
\end{proof}

We next establish an a priori estimate to the following Dirichlet
problem
\begin{equation}\label{apreq}
\begin{cases}
\mathrm{div\,}A(x,Dw)=\mathrm{div\,} \left( \frac{\varphi'(|\mathbf{f}|)}{|\mathbf{f}|} \mathbf{f} \right) & \textrm{in}\ B_{\rho} \\
w\in w_{0}+ W^{1,\varphi^{1-\delta}}_{0}( B_{\rho})
\end{cases}
\end{equation}
for the purpose of proving an existence result.

\begin{lemma}
\label{apriori} Assume \eqref{moncondi}  and \eqref{youngcondi}. Then
there exists a small constant  \linebreak $\delta_{1}=\delta_{1}(n, s_{\varphi},
\nu, L)$ such that the following holds: For any $\delta \in (0,
\delta_{1}]$, if    \linebreak $w \in
W^{1,\varphi^{1-\delta}}( B_{\rho})$ is a very weak solution to \eqref{apreq} with $w_{0}
\in W^{1,\varphi^{1-\delta}}( B_{\rho})$ and $ \mathbf{f} \in
L^{\varphi^{1-\delta}}( B_{\rho})$, then there holds
\begin{equation}\label{apresti}
 \int_{B_{\rho}} [ \varphi(|Dw|)]^{1-\delta} \, dx   \leq c \left[  \int_{B_{\rho}} [\varphi(|Dw_{0}|)]^{1-\delta} \, dx +\int_{B_{\rho}} [ \varphi(|\mathbf{f}|)]^{1-\delta} \, dx \right] ,
\end{equation}
where $c>0$ depends on  $ n,s_{\varphi}, \nu, L$.
\end{lemma}

\begin{proof}
Suppose  $\delta \in (0,\frac{1-\theta}{2} ]$ and define $v :=
w-w_{0} \in W_{0}^{1,\varphi^{1-\delta}}( B_{\rho})$. Then we apply
Lemma \ref{liptrun} to get a Lipschitz truncation $v_\lambda \in
W_{0}^{1,\infty} (B_{\rho})$ such that $v_\lambda=v$,
$Dv_\lambda=Dv$ a.e. on $B_{\rho} \backslash E_{\lambda}$, where
$$
 \quad E_{\lambda}:=  \{x \in B_{\rho} :  \left\{ \mathcal{M}([\varphi(|Dv|)]^{\theta})(x) \right\}^{\frac{1}{\theta}} >\lambda \}
$$
and we have that $\varphi(|Dv_{\lambda}|)\leq c\lambda$ on $
B_{\rho}$ for some positive constant depending on $n$ and
$s_{\varphi}$. Taking $v_{\lambda}$ as a test function in the
equation \eqref{apreq}, we have
\begin{eqnarray}\label{ii1esti1}
\notag &&  \hspace{-10mm}  \int_{B_{\rho} \backslash E_{\lambda}} \big\langle  A(x,Dw) , Dv_{\lambda}\big\rangle   \, dx = - \int_{ E_{\lambda}} \big\langle  A(x,Dw) , Dv_{\lambda}\big\rangle   \, dx \\
\notag &&  \hspace{20mm} + \int_{B_{\rho} \backslash E_{\lambda}} \big\langle  \frac{\varphi'(|\mathbf{f}|)}{|\mathbf{f}|} \mathbf{f}, Dv_{\lambda}\big\rangle   \, dx + \int_{ E_{\lambda}} \big\langle  \frac{\varphi'(|\mathbf{f}|)}{|\mathbf{f}|} \mathbf{f}, Dv_{\lambda}\big\rangle   \, dx \\
\notag &&\hspace{10mm} \leq c \left( \varphi^{-1}(\lambda)\int_{ E_{\lambda}}  \varphi ' (|Dw|)   \, dx \right. \\
 && \hspace{20mm} \left. +\int_{B_{\rho} \backslash E_{\lambda}}  \varphi ' (|\mathbf{f}|)  | Dv| \, dx +  \varphi^{-1}(\lambda)\int_{ E_{\lambda}}  \varphi ' (|\mathbf{f}|)   \, dx\right) .
\end{eqnarray}
Multiplying \eqref{ii1esti1} by $\lambda^{-(1+\delta)}$  and integrating from $0$ to $\infty$ with respect to $\lambda$, we obtain
\begin{eqnarray}\label{ii1esti2}
\notag &&  \hspace{-10mm} I_{1} : = \int_{0}^{\infty} \lambda^{-(1+\delta)} \int_{B_{\rho} \backslash E_{\lambda}} \big\langle  A(x,Dw) , Dv_{\lambda}\big\rangle   \, dx \,   d\lambda \\
\notag && \leq c_{*} \left( \int_{B_{\rho}} \left[  \int_{0}^{g(x)} \lambda^{-(1+\delta)} \varphi^{-1} (\lambda)  \varphi ' (|Dw|)    \,   d\lambda  \right]  \, dx \right.   \\
\notag && \hspace{10mm} + \int_{B_{\rho}} \left[  \int_{g(x)}^{\infty} \lambda^{-(1+\delta)}\varphi ' (|\mathbf{f}|)  | Dv|   \,   d\lambda  \right]  \, dx   \\
\notag && \hspace{15mm} \left. + \int_{B_{\rho}}  \left[  \int_{g(x)}^{\infty} \lambda^{-(1+\delta)} \varphi^{-1} (\lambda)  \varphi ' (|\mathbf{f}|)   \,   d\lambda  \right]  \, dx \right)   \\
&& \leq  c_{*} ( I_{2}+I_{3}+I_{4} ),
\end{eqnarray}
for some positive constant $c_{*}$ depending on $n,s_{\varphi}$ and
$L$, where we have defined
\begin{equation}\label{gdef2}
g(x):=  \left\{ \mathcal{M}([\varphi(|Dv|)]^{\theta})(x)
\right\}^{\frac{1}{\theta}} \quad (x \in B_{\rho}).
\end{equation}
Applying  Lemma \ref{youngineq}, $I_{2}$ can be estimated as
\begin{eqnarray}\label{ii1esti3}
\notag && \hspace{-4mm} I_{2}  = \int_{ B_{\rho}}  \varphi^{-1} (g(x)) [ g(x)]^{-\delta} \varphi ' (|Dw|)       \, dx \\
\notag && \leq c \left( \int_{ B_{\rho}}  [g(x)]^{1-\delta}   \, dx  + \int_{ B_{\rho}} [ \varphi(|Dw|)]^{1-\delta}   \, dx \right)\\
\notag && \leq c \left( \int_{B_{\rho}}  [\varphi(|Dv|)]^{1-\delta}  \, dx  + \int_{B_{\rho}}  [\varphi(|Dw|)]^{1-\delta}   \, dx \right)\\
&& \leq c \left(  \int_{ B_{\rho}} [ \varphi(|Dw|)]^{1-\delta}    \, dx  + \int_{ B_{\rho}} [ \varphi(|Dw_{0}|)]^{1-\delta}     \, dx \right) .
\end{eqnarray}
Similarly, we estimate $I_{3}$ as follows:
\begin{eqnarray}\label{ii1esti4}
\notag && \hspace{-4mm} I_{3} = \frac{1}{\delta}\int_{B_{\rho}}  [ g(x)]^{-\delta}\varphi ' (|\mathbf{f}|)  | Dv|   \, dx  \leq \frac{c}{\delta}\int_{B_{\rho}}  \varphi^{-1} (g(x))  [ g(x)]^{-\delta}\varphi ' (|\mathbf{f}|)  \\
\notag&& \leq \frac{\varepsilon}{\delta}   \int_{ B_{\rho}}  [g(x)]^{1-\delta}   \, dx   + \frac{c(\varepsilon )}{\delta}   \int_{ B_{\rho}} [ \varphi (|\mathbf{f}|)]^{1-\delta}     \, dx  \\
\notag && \leq \frac{c \varepsilon }{\delta}   \int_{ B_{\rho}}    [\varphi(|Dv|)]^{1-\delta}  \, dx   + \frac{c(\varepsilon )}{\delta}   \int_{ B_{\rho}} [ \varphi (|\mathbf{f}|)]^{1-\delta}     \, dx  \\
\notag && \leq \frac{c \varepsilon }{\delta}   \int_{ B_{\rho}}    [\varphi(|Dw|)]^{1-\delta}  \, dx   + \frac{c(\varepsilon )}{\delta} \left[   \int_{ B_{\rho}}    [\varphi(|Dw_{0}|)]^{1-\delta}  \, dx+  \int_{ B_{\rho}} [ \varphi (|\mathbf{f}|)]^{1-\delta}     \, dx \right]  .
\end{eqnarray}
In a similar manner as we did for the estimate \eqref{i1esti2}, we have
\begin{eqnarray}\label{ii1esti4}
\notag && \hspace{-12mm} I_{4} = \int_{B_{\rho}}  \left[ \int_{g(x)}^{\infty} \lambda^{-(1+\delta)} \varphi^{-1} (\lambda)  \varphi ' (|\mathbf{f}|)   \,   d\lambda  \right]   \, dx\\
\notag &&\hspace{-8mm} \leq c  \int_{0}^{\infty} \lambda^{-(1+\delta)} \varphi^{-1} (\lambda)\int_{ E_{\lambda}}  \varphi ' (|\mathbf{f}|)   \, dx    \,   d\lambda \\
\notag &&\hspace{-8mm}  \leq  c \int_{ B_{\rho}}  \varphi^{-1} (g(x)) [ g(x)]^{-\delta} \varphi ' (|\mathbf{f}|)       \, dx \\
\notag&&\hspace{-8mm}  \leq c \left(   \int_{ B_{\rho}}  [g(x)]^{1-\delta}   \, dx   +   \int_{ B_{\rho}} [ \varphi (|\mathbf{f}|)]^{1-\delta}     \, dx  \right)  \\
 &&\hspace{-8mm}  \leq c   \left( \int_{ B_{\rho}}    [\varphi(|Dw|)]^{1-\delta}  \, dx +\int_{ B_{\rho}}    [\varphi(|Dw|)]^{1-\delta}  \, dx   +    \int_{ B_{\rho}} [ \varphi (|\mathbf{f}|)]^{1-\delta}     \, dx \right) .
\end{eqnarray}
Now we estimate $I_{1}$ from below. Since $Dv_{\lambda}=Dv$ on $ B_{\rho} \backslash E_{\lambda}$, we have
\begin{eqnarray}\label{ii1esti5}
\notag &&   \hspace{-4mm} I_{1}  = \int_{ B_{\rho}}  \left[  \int_{g(x)}^{\infty} \lambda^{-(1+\delta)} \big\langle  A(x,Dw) , Dv\big\rangle     \,   d\lambda \right]   \, dx\\
\notag && = \frac{1}{\delta} \int_{B_{\rho}}   [{g(x)}]^{-\delta} \big\langle  A(x,Dw) , Dv\big\rangle     \,    dx \\
\notag && =  \frac{1}{\delta} \int_{B_{\rho}}   [{g(x)}]^{-\delta} \big\langle  A(x,Dw) - A(x,Dw_{0}) , Dv\big\rangle     \,    dx  \\
\notag && \hspace{10mm} +   \frac{1}{\delta} \int_{B_{\rho}}   [{g(x)}]^{-\delta} \big\langle  A(x,Dw_{0}) , Dv\big\rangle     \,    dx   \\
\notag && \geq \frac{\nu}{\delta} \int_{ B_{\rho}}   g^{-\delta} |V(Dw)-V(Dw_{0})|^2 \, dx  - \frac{L}{\delta} \int_{ B_{\rho}}   g^{-\delta}  \varphi' (|Dw_{0}|) |Dv|     \, dx  \\
 && = :  \frac{1}{\delta}(I_{5} - I_{6}).
\end{eqnarray}
To estimate $I_{5}$, we use Lemma \ref{maximucken} again to conclude
that ${g(x)}^{-\delta} \in  A_{1/\theta}$, and so
\begin{eqnarray}\label{ii2esti}
\notag && \hspace{-4mm} I_{5}= \nu \int_{ B_{\rho}}   g^{-\delta}  |V(Dw)-V(Dw_{0})|^2  \,    dx \\
&& \geq c \int_{ B_{\rho}}  g^{-\delta} \left\{ \mathcal{M}( |V(Dw)-V(Dw_{0})|^{2\theta}) \right\}^{\frac{1}{\theta}}      \, dx .
\end{eqnarray}
Also, we observe that for $x \in B_{\rho}$,
\begin{eqnarray}\label{gbound2}
\notag  &&  \hspace{-8mm} g(x) = \left\{ \mathcal{M}([\varphi(|Dv|)]^{\theta})(x)\right\}^{\frac{1}{\theta}} \\
 && \leq  c \left[   \left\{  \mathcal{M}( |V(Dw)-V(Dw_{0})|^{2\theta})(x) \right\}^{\frac{1}{\theta}} + \left\{  \mathcal{M}([\varphi(|Dw_{0}|)]^{\theta})(x) \right\}^{\frac{1}{\theta}} \right]  .
\end{eqnarray}
 Comparing the last two terms in \eqref{gbound2}, we have
\begin{align}\label{algebraic2}
\notag  &   \left\{  \mathcal{M}( |V(Dw)-V(Dw_{0})|^{2\theta}) \right\}^{\frac{1-\delta}{\theta}}    \\
& \leq  c   g^{-\delta} \left\{ \mathcal{M}( |V(Dw)-V(Dw_{0})|^{2\theta})\right\}^{\frac{1}{\theta}}    +  c \left\{ \mathcal{M}([\varphi(|Dw_{0}|)]^{\theta})\right\}^{\frac{1-\delta}{\theta}}
\end{align}
on $ B_{\rho}$. Therefore, combining  \eqref{ii2esti} and \eqref{algebraic2}, we have
\begin{eqnarray}\label{ii2esti2}
\notag && \hspace{-4mm} I_{5}  \geq c \int_{ B_{\rho}} \left\{  \mathcal{M}( |V(Dw)-V(Dw_{0})|^{2\theta}) \right\}^{\frac{1-\delta}{\theta}}         \, dx - c \int_{ B_{\rho}} \left\{  \mathcal{M}([\varphi(|Dw_{0}|)]^{\theta}) \right\}^{\frac{1-\delta}{\theta}}  \, dx \\
\notag &&   \geq c \int_{ B_{\rho}}  |V(Dw)-V(Dw_{0})|^{2-2\delta}        \, dx -   c \int_{ B_{\rho}} [\varphi(|Dw_{0}|)]^{1-\delta}  \, dx \\
&& \geq c \int_{ B_{\rho}}  [\varphi(|Dw|)]^{1-\delta}   \, dx - c \int_{ B_{\rho}}  [\varphi(|Dw_{0}|)]^{1-\delta}   \, dx .
\end{eqnarray}
Applying  Lemma \ref{youngineq}, we can estimates $I_6$ as follows.
\begin{eqnarray}\label{ii3esti}
\notag && \hspace{-4mm} I_{6} =  L \int_{ B_{\rho}}   [g(x)]^{-\delta}  \varphi' (|Dw_{0}|) |Dv|     \, dx \\
\notag && \leq L \int_{ B_{\rho}}   [g(x)]^{-\delta}  \varphi' (|Dw_{0}|)  \varphi^{-1} (g(x))     \, dx \\
\notag && \leq \varepsilon \int_{ B_{\rho}} [g(x)]^{1-\delta}    \, dx  + c(\varepsilon) \int_{ B_{\rho}} [ \varphi(|Dw_{0}|)]^{1-\delta}     \, dx  \\
\notag && \leq c \varepsilon \int_{ B_{\rho}} [\varphi(|Dv|)]^{1-\delta}    \, dx  + c(\varepsilon) \int_{ B_{\rho}} [ \varphi(|Dw_{0}|)]^{1-\delta}     \, dx\\
&& \leq c \varepsilon \int_{ B_{\rho}} [\varphi(|Dw|)]^{1-\delta}    \, dx  + c(\varepsilon) \int_{ B_{\rho}} [ \varphi(|Dw_{0}|)]^{1-\delta}     \, dx.
\end{eqnarray}
We finally combine \eqref{ii1esti3}, \eqref{ii2esti2}  and
\eqref{ii3esti}  to conclude
\begin{eqnarray}\label{apri1}
\notag &&\int_{ B_{\rho}}  [\varphi(|Dw|)]^{1-\delta}   \, dx  \leq c_{*} \{ \varepsilon+\delta \} \int_{ B_{\rho}} [\varphi(|Dw|)]^{1-\delta}    \, dx   \\
 && \hspace{10mm} + c_{*} \{ c(\varepsilon)+\delta \} \left[ \int_{ B_{\rho}} [ \varphi(|Dw_{0}|)]^{1-\delta}   +\int_{ B_{\rho}} [ \varphi (|\mathbf{f}|)]^{1-\delta}     \, dx  \, dx \right]
\end{eqnarray}
for some constant $c_{*}>0$ depends on $n,s_{\varphi}, \nu, L$. Taking $\varepsilon = \delta_{1}= \frac{1}{4c_{*}}$, we get the desired estimate \eqref{apresti}.
\end{proof}

Since the obtained $ W^{1,\varphi^{1-\delta}}(B_{\rho}) $-estimate
only works for a very weak solution in $
W^{1,\varphi^{1-\delta}}(B_{\rho}) $, we need to do more works in
order to get the desired gradient estimate $
W^{1,\varphi^{q}}(B_{\rho}) $ with some higher exponent $q \geq
1-\delta$.  We then prove an existence result within $
W^{1,\varphi^{1-\delta}}(B_{\rho}) $ under the condition $\delta \in
(0,\delta_{2}]$, where $\delta_{2}$ is defined as $\delta_{2}:= \min
\{\sigma, \delta_{1}\}$ so that we can use Lemma \ref{higher} and
Lemma \ref{apriori}. See also \cite[Section 7]{IS}.

\begin{corollary}\label{existence}
Assume \eqref{moncondi} and \eqref{youngcondi}. Then for all
$w_{0}\in W^{1,\varphi^{1-\delta}}  ( B_{\rho}) $ with $\delta \in
(0,\delta_{2}]$, there exists a very weak solution $w\in
W^{1,\varphi^{1-\delta}}  ( B_{\rho})$ to the problem
\begin{equation}\label{exieq}
\begin{cases}
\mathrm{div\,}A(x,Dw)=0 & \textrm{in}\ B_{\rho} \\
w\in w_{0}+ W^{1,\varphi^{1-\delta}}_{0}( B_{\rho}),
\end{cases}
\end{equation}
such that we have the estimates
\begin{equation}\label{exiesti}
\int_{B_{\rho}}  [\varphi(|Dw|)]^{1-\delta} \, dx   \leq c   \int_{B_{\rho}} [\varphi(|Dw_{0}|)]^{1-\delta} \, dx
\end{equation}
where the positive constant $c$ depends on  $ n, s_{\varphi}, \nu $ and $ L$.
\end{corollary}

\begin{proof}
Let  $ \delta \in (0, \delta_{2}]$ and take a sequence of functions $w_{0,k} \in C^{\infty} ( \overline{B_{\rho}})$ such that
\begin{equation}\label{w0kconv}
w_{0,k} \rightarrow w_{0} \quad \textrm{strongly in} \ W^{1,\varphi^{1-\delta}}  ( B_{\rho}).
\end{equation}
Then, there exists a unique weak solution $w_{k} \in W^{1,\varphi}  ( B_{\rho})$ to the boundary value problem
\begin{equation}\label{seqeq}
\begin{cases}
\mathrm{div\,}A(x,Dw_{k})=0 & \textrm{in}\ B_{\rho} \\
w_{k}\in w_{0,k}+ W^{1,\varphi}_{0}( B_{\rho}).
\end{cases}
\end{equation}
 The solvability follows from \cite[Section 2]{Sh}. According to Lemma \ref{apriori}, we get a uniform bound
\begin{equation}\label{seqbound}
\int_{B_{\rho}}  [\varphi(|Dw_{k}|)]^{1-\delta} \, dx   \leq c   \int_{B_{\rho}} [\varphi(|Dw_{0,k}|)]^{1-\delta} \, dx  \leq c   \int_{B_{\rho}} [\varphi(|Dw_{0}|)]^{1-\delta} \, dx,
\end{equation}
for sufficiently large $k$. Combining \eqref{seqbound} with Lemma \ref{poincare}, we conclude that the norm  $\|w_{k}\|_{W^{1,\varphi^{1-\delta}}  (B_{\rho})}$ is uniformly bounded. Thus, passing to a subsequence, we have a function $w \in {W^{1,\varphi^{1-\delta}} }( B_{\rho})$ such that
\begin{equation}\label{wkconv}
w_{k} \rightharpoonup w \quad \textrm{weakly in} \ W^{1,\varphi^{1-\delta}}  ( B_{\rho}),  \quad w_{k} \rightarrow w \quad \textrm{strongly in} \ L^{\varphi^{1-\delta}} (B_{\rho}).
\end{equation}
Similarly, we consider a sequence $\tilde{w}_{k} := w_{k}-w_{0,k}$ and take $k \rightarrow \infty $ to observe that $w\in w_{0}+ W^{1,\varphi^{1-\delta}}_{0}( B_{\rho})$.

Now, fix any test function $\phi \in C_{0}^{\infty}( B_{\rho})$ and choose an open set $U$ such that $\textrm{supp} \ \phi \subset U \subset \subset  B_{\rho} $. Then in light of Lemma \ref{higher} we have a uniform bound of $\|w_{k}\|_{W^{1,\varphi}  ( U)}$.  Passing to a subsequence, we find that
\begin{equation}\label{wkconv2}
w_{k} \rightharpoonup w \quad \textrm{weakly in} \ W^{1,\varphi} (U), \quad w_{k} \rightarrow w \quad \textrm{strongly in} \ L^{\varphi} (U).
\end{equation}
We recall the Minty-Browder technique to conclude that $w$ is a weak solution to $\mathrm{div\,}A(x,Dw)=0$ in $U$, which means that
\begin{eqnarray}\label{solutionpf}
\notag && \hspace{-3mm} 0 = \int_{U} \big\langle  A(x,Dw) , D\phi \big\rangle   \, dx =   \int_{B_{\rho}} \big\langle  A(x,Dw) , D\phi \big\rangle   \, dx.
\end{eqnarray}
Since $\phi$ is arbitrary, $w$ becomes a very weak solution of \eqref{exieq}.
\end{proof}

\section{Proof of Theorem \ref{mainthm}}\label{Sec5}
This section is devoted to proving the main result, Theorem
\ref{mainthm}. To do this, we need to compare a solution $u$ to
\eqref{maineq} under consideration with a solution $w$ to the
following homogeneous problem
\begin{equation}\label{refeq}
\begin{cases}
\mathrm{div\,}A(x,Dw)=0 & \textrm{in}\ B_{2r} \\
w\in u+ W^{1,\varphi^{1-\delta_{0}}}_{0}( B_{2r}),
\end{cases}
\end{equation}
where $\delta_{0}>0$ is to be selected later. We further assume that
\begin{equation}\label{intlambda}
\Xint-_{B_{2r}} [\varphi(|Du|)]^{1-\delta_{0}} \, dx \leq  \Lambda,  \quad   \hspace{0.25em} \Xint-_{B_{2r}}  [ \varphi(|\mathbf{f}|)]^{1-\delta_{0}} \, dx \leq \delta_{0} \Lambda
\end{equation}
for some $\Lambda>0$.
These assumptions will be made during an exit time argument below. We will frequently use the universal constants $\sigma, \delta_{1}$ and $\delta_{2}$ given in the previous section.

\begin{lemma}\label{compari}
Assume  \eqref{moncondi} and \eqref{youngcondi}. Then for any
$\varepsilon \in (0,1]$, there exists a constant $ \delta_{0} = \delta_{0}(n,
s_{\varphi}, \nu, L, \varepsilon ) $ $ \in {(0, \frac{\sigma}{2}]}$
such that the following holds: For any very weak solution $u \in
W^{1,\varphi^{1-\delta_{0}}} (\Omega)$ to \eqref{maineq} with
$\mathbf{f} \in W^{1,\varphi^{1-\delta_{0}}}(\Omega)$, if
\eqref{intlambda} holds, then there exists a very weak solution $w
\in W^{1,\varphi^{1-\delta_{0}}}  ( B_{2r}  ) $ to the equation
\eqref{refeq} such that
\begin{align}\label{compa}
 \Xint-_{B_{2r}}|V(Du)-V(Dw)|^{2-2\delta_{0}} \, dx \leq \varepsilon \Lambda
\end{align}
with the estimate
\begin{align}\label{compenergy}
 \Xint-_{B_{2r}} [\varphi(|Dw|)]^{1-\delta_{0}} \, dx  \leq  c \Lambda,  \quad   \left(  \hspace{0.25em} \Xint-_{B_{r}} [\varphi(|Dw|)]^{1+\sigma} \, dx  \right)^{\frac{1-\delta_0}{1+\sigma}}  \leq  c \Lambda,
\end{align}
where the constant $c>0$ depends only on  $ n, s_{\varphi}, \nu$ and $L$.
\end{lemma}

\begin{proof}
Let $0 <\delta_{0} \leq \frac{\delta_{2}}{2} $. Then according to Lemma \ref{higher} and Corollary \ref{existence}, there exists a very weak solution $w \in W^{1,\delta_{0}}(B_{2r})$ to \eqref{refeq} such that the estimates \eqref{compenergy} holds, since we have
\begin{eqnarray}\label{holder}
\notag && \hspace{-10mm}  \left(  \hspace{0.25em} \Xint-_{B_{r}} [\varphi(|Dw|)]^{1+\sigma} \, dx  \right)^{\frac{1-\delta_0}{1+\sigma}} \leq   c^{1-\delta_{0}} \left(  \hspace{0.25em} \Xint-_{B_{2r}} [\varphi(|Dw|)]^{1-\sigma} \, dx  \right)^{\frac{1-\delta_0}{1-\sigma}} \\
 && \hspace{35mm}\leq  c  \hspace{0.25em}  \Xint-_{B_{2r}} [\varphi(|Dw|)]^{1-\delta_{0}} \, dx   \leq  c\Lambda.
\end{eqnarray}

We next write $v := u-w \in W_{0}^{1,\varphi^{1- \delta_{0}}}(
B_{2r} )$ and take a Lipschitz truncation $v_\lambda \in
W_{0}^{1,\infty} ( B_{2r} )$ with $\lambda>0$  in light of Lemma
\ref{liptrun}. We have $v_\lambda=v$, $Dv_\lambda=Dv$ a.e. on
$B_{2r}\backslash E_{\lambda}$ for
$$
 \quad E_{\lambda}:=  \{x \in B_{2r}  :  \left\{ \mathcal{M}([\varphi(|Dv|)]^{\theta})(x) \right\}^{\frac{1}{\theta}} >\lambda
 \},
$$
and $\varphi(|Dv_{\lambda}|)\leq c\lambda$ for $x \in B_{2r} $.  Taking $v_{\lambda}$  as a test function to both   \eqref{maineq} and \eqref{refeq}, we get \pagebreak
\begin{eqnarray}\label{iii1esti1}
\notag &&  \hspace{-10mm}  \int_{B_{2r}  \backslash E_{\lambda}} \big\langle  A(x,Du)-A(x,Dw) , Dv_{\lambda}\big\rangle   \, dx \\
\notag&&\hspace{-5mm}= - \int_{ E_{\lambda}} \big\langle  A(x,Du)-A(x,Dw) , Dv_{\lambda}\big\rangle   \, dx + \int_{ B_{2r} } \big\langle  \frac{\varphi'(|\mathbf{f}|)}{|\mathbf{f}|} \mathbf{f}, Dv_{\lambda}\big\rangle   \, dx \\
\notag &&\hspace{-10mm} \leq c \left( \varphi^{-1}(\lambda)\int_{ E_{\lambda}}  \varphi ' (|Du|)   \, dx + \varphi^{-1}(\lambda)\int_{ E_{\lambda}}  \varphi ' (|Dw|)   \, dx \right. \\
 && \hspace{10mm} \left. +\int_{B_{2r} \backslash E_{\lambda}}  \varphi ' (|\mathbf{f}|)  | Dv| \, dx +  \varphi^{-1}(\lambda)\int_{ E_{\lambda}}  \varphi ' (|\mathbf{f}|)   \, dx\right).
\end{eqnarray}
Multiplying \eqref{iii1esti1} by $\lambda^{-(1+ \delta_{0})}$  and integrating from $0$ to $\infty$ with respect to $\lambda$, along with the similar calculations as in \eqref{ii1esti1}, we discover that
\begin{eqnarray}\label{iii1esti2}
\notag &&  \hspace{-5mm}  \delta_{0} I_{1} : = \delta_{0} \int_{0}^{\infty} \lambda^{-(1+ \delta_{0})} \int_{B_{2r} \backslash E_{\lambda}} \big\langle  A(x,Du)-A(x,Dw) , Dv_{\lambda}\big\rangle   \, dx \,   d\lambda \\
\notag &&  \hspace{3mm} \leq c \{ \varepsilon_{1} +\delta_{0} \}   \left[    \int_{ B_{2r} }    [\varphi(|Du|)]^{1-\delta_{0}}  \, dx  + \int_{ B_{2r} }    [\varphi(|Dw|)]^{1-\delta_{0}}  \, dx  \right] \\
&&\hspace{15mm} + c \{ c( \varepsilon_{1}) +\delta_{0} \}   \int_{ B_{2r} } [ \varphi (|\mathbf{f}|)]^{1-\delta_{0}}     \, dx   \\
\notag &&  \hspace{3mm} \leq c \{ \varepsilon_{1}+ c( \varepsilon_{1})\delta_{0} \} \Lambda r^n.
\end{eqnarray}
We first find a lower bound of $I_{1}$. Defining ${g(x)}$ as
\begin{equation}\label{gdef3}
g(x):=  \left\{ \mathcal{M}([\varphi(|Dv|)]^{\theta})(x) \right\}^{\frac{1}{\theta}} \quad (x \in B_{2r}),
\end{equation}
we again have $[{g(x)}]^{-\delta_{0}} \in  A_{1/\theta}$ by Lemma \ref{maximucken}. Then it follows that
\begin{eqnarray}\label{iii1esti3}
\notag &&   \hspace{-4mm} I_{1}  = \int_{ B_{2r} }  \int_{g(x)}^{\infty} \lambda^{-(1+\delta_{0})} \big\langle  A(x,Du)-A(x,Dw) , Dv\big\rangle     \,   d\lambda   \, dx\\
\notag && = \frac{1}{\delta_{0}} \int_{B_{2r} }   [{g(x)}]^{-\delta_{0}} \big\langle  A(x,Du)-A(x,Dw) , Dv\big\rangle     \,    dx \\
\notag  && \geq \frac{\nu}{\delta_{0}} \int_{ B_{2r} }   g^{-\delta_{0}} |V(Du)-V(Dw)|^2 \, dx \\
&& \geq \frac{c}{\delta_{0}} \int_{ B_{2r} }  g^{-\delta_{0}} \left\{ \mathcal{M}( |V(Du)-V(Dw)|^{2\theta})\right\}^{\frac{1}{\theta}}      \, dx.
\end{eqnarray}
We now observe from \eqref{phitvt} that  for any $x \in B_{2r} $ and for any $\tilde{\varepsilon} \in (0,1]$,
\begin{eqnarray}\label{gbound3}
\notag  &&  \hspace{-10mm} g(x) =  \left\{ \mathcal{M}([\varphi(|Dv|)]^{\theta}) (x) \right\}^{\frac{1}{\theta}}\\
\notag   && \hspace{-2mm} \leq \left\{  \mathcal{M}( \left[ c( \tilde{\varepsilon} ) |V(Du)-V(Dw)|^{2\theta}+ \tilde{\varepsilon} ^{\theta}[\varphi(|Du|)]^{\theta} \right] )(x) \right\}^{\frac{1}{\theta}} \\
 && \hspace{-2mm} \leq  c ( \tilde{\varepsilon} )  \left\{ \mathcal{M}( |V(Du)-V(Dw)|^{2\theta})(x)\right\}^{\frac{1}{\theta}}  +  c \tilde{\varepsilon} \left\{ \mathcal{M}([\varphi(|Du|)]^{\theta})(x)\right\}^{\frac{1}{\theta}} .
\end{eqnarray}
  Consequently, it follows that
\begin{align}\label{algebraic3}
\notag  &  \left\{ \mathcal{M}( |V(Du)-V(Dw)|^{2\theta})\right\}^{\frac{1-\delta_{0}}{\theta}}    \\
& \hspace{5mm} \leq  c( \varepsilon_{2} )  g^{-\delta_{0}} \left\{ \mathcal{M}( |V(Du)-V(Dw)|^{2\theta}) \right\}^{\frac{1}{\theta}}    +  \varepsilon_{2} \left\{ \mathcal{M}([\varphi(|Du|)]^{\theta}) \right\}^{\frac{1-\delta_{0}}{\theta}}
\end{align}
for any $\varepsilon_{2} \in (0,1]$. Recalling \eqref{iii1esti3}, we have
\begin{eqnarray}\label{iii2esti2}
\notag && \hspace{-4mm} I_{1}  \geq \frac{1}{ c(\varepsilon_{2})}\int_{B_{2r} }  \left\{ \mathcal{M}( |V(Du)-V(Dw)|^{2\delta_{0}}) \right\}^{\frac{1-\delta}{\theta}}         \, dx \\
\notag && \hspace{10mm}  - \frac{  \varepsilon_{2} }{ c(\varepsilon_{2})}   \int_{ B_{2r} } \left\{ \mathcal{M}([\varphi(|Du|)]^{\theta}) \right\}^{\frac{1-\delta_{0}}{\theta}}  \, dx \\
\notag  &&   \geq \frac{1}{ c(\varepsilon_{2})} \int_{ B_{2r} }  |V(Du)-V(Dw)|^{2-2\delta_{0}}        \, dx -   \frac{   \varepsilon_{2} }{ c(\varepsilon_{2})}  \int_{ B_{2r} } [\varphi(|Du|)]^{1-\delta_{0}}  \, dx \\
 &&   \geq \frac{1}{ c(\varepsilon_{2})} \int_{ B_{2r} }  |V(Du)-V(Dw)|^{2-2\delta_{0}}        \, dx -   \frac{  \varepsilon_{2} }{ c(\varepsilon_{2})}  \Lambda  .
\end{eqnarray}
Combining \eqref{intlambda}, \eqref{iii1esti3} and \eqref{iii2esti2}, we conclude that
\begin{eqnarray}\label{compenergy2}
\notag && \int_{ B_{2r} }  |V(Du)-V(Dw)|^{2-2\delta_{0}}     \, dx \leq  (   \varepsilon_{2}+ c (\varepsilon_{2})  \varepsilon_{1} +  c (\varepsilon_{1}) c( \varepsilon_{2})\delta_{0} )\Lambda   \leq 3\varepsilon \Lambda,
\end{eqnarray}
by taking $\varepsilon_{2} = \varepsilon$, $\varepsilon_{1}
=\frac{\varepsilon}{c (\varepsilon_{2})} $ and then $ \delta_{0}
:=\min{\{ \frac{\varepsilon}{c (\varepsilon_{1})c (\varepsilon_{2})}
, \frac{\delta_{2}}{2}\}}$ for the last inequality. This completes
the proof.
\end{proof}

We are in a position to prove the main result of the paper.

\begin{proof}[Proof of Theorem \ref{mainthm}] Our proof is based on the harmonic
analysis-free technique which was first introduced in \cite{AM}. We
are under the same assumption as in Theorem \ref{mainthm}, where
$\delta_{0}$ is given in Lemma \ref{compari}. Notice that if we
choose $\varepsilon$ in Lemma \ref{compari} depending only on
$n,s_{\varphi},\nu$ and $L$, accordingly $\delta_{0}$ is to be
selected depending only on $n,s_{\varphi},\nu$ and $L$. First we fix
any ball $B_{R}(z) \subset \mr^n$ and write upper level sets as
$$
E^{s}_{\Lambda} := \{x \in B_{s}(z) :  [\varphi(|Du|)]^{1-\delta_{0}}>\Lambda \}  \ \ \ \left(\Lambda>0\right),
$$
for $ \frac{R}{2} \leq s \leq R$.  Consider the concentric balls
$B_{r_1}(z)$ and $ B_{r_2}(z)$ with $ \frac{R}{2} \leq r_{1} < r_{2}
\leq R$. Then for each $y\in E^{r_{1} }_{\Lambda} $, we define a
continuous function $\Phi_{y} : (0,r_{2}-r_{1}] \rightarrow
[0,\infty)$ by
\begin{equation}\label{intaverge}
\Phi_{y}(\rho):=   \hspace{0.25em}  \Xint-_{B_{\rho}(y)} \left\{ [\varphi(|Du|)]^{1-\delta_{0}}
+ \frac{ [ \varphi(|\mathbf{f}|)]^{1-\delta_{0}}}{\delta_{0}} \right\} \, dx.
\end{equation}
Applying Lebesgue differentiation theorem, we get
\begin{equation}\label{almostlambda}
\lim_{\rho \rightarrow 0}\Phi_{y}(\rho)=  [\varphi(|Du|)]^{1-\delta_{0}}
+ \frac{ [ \varphi(|\mathbf{f}|)]^{1-\delta_{0}}}{\delta_{0}} \geq [\varphi(|Du|)]^{1-\delta_{0}} > \Lambda
\end{equation}
for a.e. $y \in E^{r_{1} }_{\Lambda}$. Note that if $ \frac{r_{2}-r_{1} }{10} \leq \rho \leq r_{2}-r_{1}  $, then
\begin{equation}\label{lambdazero}
 \Phi_{y}(\rho) \leq  \frac{10^n r_{2}^n}{(r_{2}-r_{1})^{n}} \hspace{0.25em} \Xint-_{B_{r_{2}}(z)}  \left\{ [\varphi(|Du|)]^{1-\delta_{0}}
+ \frac{ [ \varphi(|\mathbf{f}|)]^{1-\delta_{0}}}{\delta_{0}} \right\} \, dx =: \Lambda_{0}.
\end{equation}
For $\Lambda> \Lambda_{0}$, since $\Phi_{y}$ is continuous and
$\lim_{\rho \rightarrow 0}\Phi_{y}(\rho) =\Lambda > \Lambda_{0}$,
there exists an exit time radius $\rho_{y} \in (0,  \frac{
r_{2}-r_{1} }{10} )$ such that
\begin{equation}\label{exittime}
\Phi_{y}(\rho_{y}) = \Lambda \quad \textrm{and} \quad \Phi_{y}(\rho)<\Lambda \quad \textrm{if} \enspace  \rho \in (\rho_{y},  r_{2}-r_{1} ].
\end{equation}
Now we consider the family $\{B_{\rho_{y}}(y) : y \in
E^{r_{1} }_{\Lambda}\}$ which covers the set $E^{r_{1} }_{\Lambda}$. By Vitali's
covering lemma, we find a countable family of disjoint sets $\{B_{\rho_{i}}(y_{i}) : y_{i} \in E^{r_{1} }_{\Lambda} \}$ such that
\begin{equation}\label{covering}
E^{r_{1} }_{\Lambda} \subset \bigcup_{i\geq 1} B_{5\rho_{i}}(y_{i}) \cup \textrm{negligible set} ,
\end{equation}
where we have denoted $\rho_{i}=\rho_{y_{i}}$. We write $B^{i}_{r}:= B_{5\rho_{i}}(y_{i})$ and so $B^{i}_{2r} = B_{10\rho_{i}}(y_{i})$. It follows from \eqref{exittime} that
\begin{equation}\label{bdylambda}
\Xint-_{B^{i}_{2r}}  [\varphi(|Du|)]^{1-\delta_{0}} \, dx \leq \Lambda,  \quad  \hspace{0.25em} \Xint-_{B^{i}_{2r}}  [ \varphi(|\mathbf{f}|)]^{1-\delta_{0}} \, dx \leq \delta_{0}\Lambda,
\end{equation}
which verifies the assumptions \eqref{intlambda}. According to Lemma \ref{compari}, we find that for $\Lambda> \Lambda_{0}$ and for any $T \geq 1$,
\begin{eqnarray}\label{loclevelsetesti1}
\nonumber &&  \hspace{-5mm} \int_{\{x \in B^{i}_{r} : [\varphi(|Du|)]^{1-\delta_{0}}>T \Lambda \}}  [\varphi(|Du|)]^{1-\delta_{0}}  \, dx \\
\nonumber &&   \hspace{5mm} \leq   2 \int_{\{x \in B^{i}_{r} : [\varphi(|Du|)]^{1-\delta_{0}}>T\Lambda \}}  |V(Du)-V(Dw)|^{2-2\delta_{0}} \, dx   \\
 &&  \hspace{15mm}  +2 \int_{\{x \in B^{i}_{r} : [\varphi(|Du|)]^{1-\delta_{0}}>T\Lambda \}}  [\varphi(|Dw|)]^{1-\delta_{0}}  \, dx ,
\end{eqnarray}
where the last term can be estimated as
\begin{eqnarray}\label{loclevelsetesti1}
\nonumber &&\hspace{-10mm} \int_{\{x \in B^{i}_{r} : [\varphi(|Du|)]^{1-\delta_{0}}>T\Lambda \}}  [\varphi(|Dw|)]^{1-\delta_{0}}  \, dx  \\
\nonumber &&  \hspace{-5mm} \leq |\{x \in B^{i}_{r} : [\varphi(|Du|)]^{1-\delta_{0}}>T\Lambda \}|^{\frac{\delta_{0}+\sigma}{1+\sigma}}  \left( \int_{B^{i}_{r}} [\varphi(|Dw|)]^{1+\sigma} \, dx  \right)^{\frac{1-\delta_0}{1+\sigma}} \\
\nonumber &&   \hspace{-5mm} \leq  c  \left( \frac{1}{T\Lambda}\int_{\{x \in B^{i}_{r} : [\varphi(|Du|)]^{1-\delta_{0}}>T \Lambda \}}  [\varphi(|Du|)]^{1-\delta_{0}}  \, dx \right)^{\frac{\delta_{0}+\sigma}{1+\sigma}}   \Lambda  |B^{i}_{2r}|^{\frac{1-\delta_{0}}{1+\sigma}} \\
\nonumber &&   \hspace{-5mm} =  c  \left( \int_{\{x \in B^{i}_{r} : [\varphi(|Du|)]^{1-\delta_{0}}>T \Lambda \}}  [\varphi(|Du|)]^{1-\delta_{0}}  \, dx \right)^{\frac{\delta_{0}+\sigma}{1+\sigma}}  \left[  T^{-\frac{\delta_{0}+\sigma}{1-\delta_{0}}} \Lambda  |B^{i}_{2r}| \right]^{\frac{1-\delta_{0}}{1+\sigma}} \\
 &&   \hspace{-5mm}  \leq   \frac{1}{2} \int_{\{x \in B^{i}_{r} : [\varphi(|Du|)]^{1-\delta_{0}}>T \Lambda \}}  [\varphi(|Du|)]^{1-\delta_{0}}  \, dx +  c T^{-\frac{\delta_{0}+\sigma}{1-\delta_{0}}}\Lambda |B^{i}_{2r}|,
\end{eqnarray}
where we have used Chebyshev's inequality for the second inequality. Then we find
\begin{eqnarray}\label{loclevelsetesti3}
\nonumber &&  \hspace{-15mm}  \int_{\{x \in B^{i}_{r} : [\varphi(|Du|)]^{1-\delta_{0}}>T \Lambda \}}  [\varphi(|Du|)]^{1-\delta_{0}}  \, dx \leq ( 4  \varepsilon + c^{*} T^{-\frac{\delta_{0}+\sigma}{1-\delta_{0}}} ) \Lambda |B^{i}_{2r}|  \\
&& \hspace{35mm}  \leq  40^{n}( \varepsilon + c^{*}  T^{-\frac{\delta_{0}+\sigma}{1-\delta_{0}}} ) \Lambda |B_{\rho_{i}}(y_{i}) |
\end{eqnarray}
 for some $c^{*}$ depending only on $n, s_{\varphi}, \nu, L$. Now we estimate
\begin{eqnarray}\label{loclevelsetesti4}
\nonumber &&  \hspace{-15mm} |B_{\rho_{i}}(y_{i}) | =  \frac{1}{\Lambda}  \int_{B_{\rho_{i}}(y_{i}) } \left\{ [\varphi(|Du|)]^{1-\delta_{0}}
+ \frac{ [ \varphi(|\mathbf{f}|)]^{1-\delta_{0}}}{\delta_{0}} \right\} \, dx \\
\nonumber &&    \hspace{-10mm}  \leq    \frac{1}{\Lambda}  \int_{\{x \in B_{\rho_{i}}(y_{i}) : [\varphi(|Du|)]^{1-\delta_{0}}> \frac{\Lambda}{4} \} }  [\varphi(|Du|)]^{1-\delta_{0}}  \, dx \\
&& + \frac{1}{\Lambda}  \int_{\{x \in B_{\rho_{i}}(y_{i}) : [\varphi(|\mathbf{f}|)]^{1-\delta_{0}}> \frac{\delta_{0}\Lambda}{4} \} } \frac{ [ \varphi(|\mathbf{f}|)]^{1-\delta_{0}}}{\delta_{0}}   \, dx + \frac{ |B_{\rho_{i}}(y_{i}) |}{2}.
\end{eqnarray}
Combining \eqref{loclevelsetesti3} with \eqref{loclevelsetesti4}, we find \pagebreak
\begin{eqnarray}\label{loclevelsetesti5}
\nonumber &&   \hspace{-5mm}   \int_{\{x \in B^{i}_{r} : [\varphi(|Du|)]^{1-\delta_{0}}>T \Lambda \}}  [\varphi(|Du|)]^{1-\delta_{0}}  \, dx \\
\nonumber &&    \leq   80^{n}( \varepsilon + c^{*} T^{-\frac{\delta_{0}+\sigma}{1-\delta_{0}}} )  \left[  \int_{\{x \in B_{\rho_{i}}(y_{i}) : [\varphi(|Du|)]^{1-\delta_{0}}> \frac{\Lambda}{4} \} }  [\varphi(|Du|)]^{1-\delta_{0}}  \, dx \right. \\
&&  \hspace{5mm}   \left. +   \int_{\{x \in B_{\rho_{i}}(y_{i}) : [\varphi(|\mathbf{f}|)]^{1-\delta_{0}}> \frac{\delta_{0}\Lambda}{4} \} }  \frac{ [ \varphi(|\mathbf{f}|)]^{1-\delta_{0}}}{\delta_{0}}    \, dx \right] .
\end{eqnarray}
Write upper level sets of $[\varphi(|\mathbf{f}|)]^{1-\delta_{0}}$ by
$$
\mathcal{E}^{s}_{\Lambda} = \{x \in  B_{s}(z) : [\varphi(|\mathbf{f}|)]^{1-\delta_{0}}>\Lambda \}.
$$
Recall that the set $E^{r_{1}}_{T\Lambda}\subset E^{r_{1}}_{\Lambda}$ is covered by the family $\{ B^{i}_{r} \}_{i=1}^{\infty} $, where $\{ B_{\rho_{i}}(y_{i}) \}_{i=1}^{\infty}$ is a disjoint family of balls. Summing up over the covering $\{ B^{i}_{r} \}_{i=1}^{\infty}$, we get
\begin{eqnarray}\label{levelsetesti1}
\nonumber && \hspace{-10mm}   \int_{E^{r_{1}}_{T\Lambda}}  [\varphi(|Du|)]^{1-\delta_{0}} \, dx  \\
\nonumber && \hspace{-5mm}    \leq  \sum\limits_{i \geq 1}  \int_{\{x \in B^{i}_{r} : [\varphi(|Du|)]^{1-\delta_{0}}>T \Lambda \}}  [\varphi(|Du|)]^{1-\delta_{0}} \, dx  \\
\nonumber &&    \hspace{-5mm}  \leq      80^{n}( \varepsilon + c^{*} T^{-\frac{\delta_{0}+\sigma}{1-\delta_{0}}} )    \left[  \sum\limits_{i \geq 1}  \int_{\{x \in B_{\rho_{i}}(y_{i}) : [\varphi(|Du|)]^{1-\delta_{0}}> \frac{\Lambda}{4} \} }  [\varphi(|Du|)]^{1-\delta_{0}}  \, dx \right. \\
\nonumber &&  \hspace{25mm} \left. +  \sum\limits_{i \geq 1}   \int_{\{x \in B_{\rho_{i}}(y_{i}) : [\varphi(|\mathbf{f}|)]^{1-\delta_{0}}> \frac{\delta_{0}\Lambda}{4} \} }  \frac{ [ \varphi(|\mathbf{f}|)]^{1-\delta_{0}}}{\delta_{0}}    \, dx \right] \\
\nonumber &&    \hspace{-5mm}  \leq   80^{n}( \varepsilon + c^{*} T^{-\frac{\delta_{0}+\sigma}{1-\delta_{0}}} )   \left[   \int_{ \bigcup_{i\geq 1} \{x \in B_{\rho_{i}}(y_{i}) : [\varphi(|Du|)]^{1-\delta_{0}}> \frac{\Lambda}{4} \} }  [\varphi(|Du|)]^{1-\delta_{0}}  \, dx \right. \\
\nonumber &&  \hspace{25mm} \left. +     \int_{ \bigcup_{i\geq 1} \{x \in B_{\rho_{i}}(y_{i}) : [\varphi(|\mathbf{f}|)]^{1-\delta_{0}}> \frac{\delta_{0}\Lambda}{4} \} }  \frac{ [ \varphi(|\mathbf{f}|)]^{1-\delta_{0}}}{\delta_{0}}    \, dx \right] \\
\nonumber &&      \hspace{-5mm}     \leq  80^{n}( \varepsilon + c^{*} T^{-\frac{\delta_{0}+\sigma}{1-\delta_{0}}} )   \left[ \int_{E^{r_{2}}_{\Lambda/4}} [\varphi(|Du|)]^{1-\delta_{0}} \, dx +   \int_{\mathcal{E}^{r_{2}}_{\delta_{0}\Lambda/4} } \frac{ [ \varphi(|\mathbf{f}|)]^{1-\delta_{0}}}{\delta_{0}}   \, dx \right].
\end{eqnarray}
Thus, change of variable with respect to $\Lambda$ leads to
\begin{align}\label{levelsetesti2}
\nonumber  &  \hspace{-3mm}   \int_{E^{r_{1}}_{\Lambda}} [\varphi(|Du|)]^{1-\delta_{0}}  \, dx     \\
  &   \leq  80^{n}( \varepsilon + c^{*} T^{-\frac{\delta_{0}+\sigma}{1-\delta_{0}}} )   \left[ \int_{E^{r_{2}}_{\Lambda/4T}} [\varphi(|Du|)]^{1-\delta_{0}} \, dx +   \int_{\mathcal{E}^{r_{2}}_{\delta_{0}\Lambda/4T} } \frac{ [ \varphi(|\mathbf{f}|)]^{1-\delta_{0}}}{\delta_{0}}   \, dx \right]
\end{align}
for $\Lambda> T\Lambda_0$.

 \pagebreak

We next introduce truncation functions
$$
\left[ [\varphi(|Du|)]^{1-\delta_{0}}\right]_{t}:=\min \{ [\varphi(|Du|)]^{1-\delta_{0}} , t \} \quad (t>0) .
$$
For $t > T\Lambda_0 =: t_{0}$ and  $q \in [1-\delta_{0}, 1+\delta_{0}]$, Fubini's theorem gives
\begin{eqnarray}\label{fubiniinteg1}
\nonumber &&   \hspace{-5mm} \int_{B_{r_{1}}(z)} [\varphi(|Du|)]^{1-\delta_{0}}  \left[ [\varphi(|Du|)]^{1-\delta_{0}}\right]_{t}^{\frac{q-1+\delta_{0}}{1-\delta_{0}}} \, dx  \\
\nonumber &&= \left( \frac{q-1+\delta_{0} }{1-\delta_{0}} \right)  \int_0^{t} \Lambda^{\frac{q}{1-\delta_{0}}-2} \int_{E^{r_{1}}_{\Lambda}} [\varphi(|Du|)]^{1-\delta_{0}} \, dx \, d\Lambda \\
\nonumber && \leq \left( \frac{q-1+\delta_{0} }{1-\delta_{0}} \right)\int_0^{t_{0}} \Lambda^{\frac{q}{1-\delta_{0}}-2} \int_{E^{r_{1}}_{\Lambda}} [\varphi(|Du|)]^{1-\delta_{0}} \, dx  \, d\Lambda \\
\nonumber && \hspace{15mm} + \left( \frac{q-1+\delta_{0} }{1-\delta_{0}} \right) \int_{t_{0}}^{t} \Lambda^{\frac{q}{1-\delta_{0}}-2} \int_{E^{r_{1}}_{\Lambda}} [\varphi(|Du|)]^{1-\delta_{0}} \, dx \, d\Lambda  \\
\nonumber && \leq (T\Lambda_0)^{\frac{q}{1-\delta_{0}}-1}  \int_{B_{r_{2}}(z)} [\varphi(|Du|)]^{1-\delta_{0}} \, dx  \\
\nonumber && \hspace{15mm} + \left( \frac{q-1+\delta_{0} }{1-\delta_{0}} \right) \int_{t_{0}}^{t} \Lambda^{\frac{q-2+2\delta_{0}}{1-\delta_{0}}} \int_{E^{r_{1}}_{\Lambda}} [\varphi(|Du|)]^{1-\delta_{0}} \, dx \, d\Lambda  \\
\nonumber &&  \hspace{-2.5mm}  \underset{\eqref{lambdazero}}{\leq}  T^{\frac{q}{1-\delta_{0}}-1} \Lambda_{0}^{\frac{q}{1-\delta_{0}}}  |B_{r_{2}}(z)| \\
\nonumber && \hspace{-5mm} + 80^{n}( \varepsilon + c^{*} T^{-\frac{\delta_{0}+\sigma}{1-\delta_{0}}} )  \left[ \left( \frac{q-1+\delta_{0} }{1-\delta_{0}} \right) \int_{t_0}^{t}  \Lambda^{\frac{q}{1-\delta_{0}}-2}  \int_{E^{r_{2}}_{\Lambda/4T}} [\varphi(|Du|)]^{1-\delta_{0}} \, dx \, d\Lambda \right.\\
&& \hspace{15mm} \left. +  \left( \frac{q-1+\delta_{0} }{1-\delta_{0}} \right) \int_{t_0}^{t}  \Lambda^{\frac{q}{1-\delta_{0}}-2}  \int_{\mathcal{E}^{r_{2}}_{\delta_{0}\Lambda/4T} } [ \varphi(|\mathbf{f}|)]^{1-\delta_{0}}   \, dx  \, d\Lambda \right] .
\end{eqnarray}
To estimate the last two integrals, we again use Fubini's theorem
along with change of variable, to discover that
\begin{eqnarray}\label{fubiniinteg2}
\nonumber &&   \hspace{-15mm}   \int_{t_0}^{t}  \Lambda^{\frac{q}{1-\delta_{0}}-2}   \int_{E^{r_{2}}_{\Lambda/4T}} [\varphi(|Du|)]^{1-\delta_{0}} \, dx \, d\Lambda  \\
\nonumber &&  \hspace{-10mm} \leq (4T)^{\frac{q}{1-\delta_{0}}-1}  \left( \frac{1-\delta_{0}}{q-1+\delta_{0} } \right)  \int_{B_{r_{2}}(z)} [\varphi(|Du|)]^{1-\delta_{0}}  [\varphi(|Du|)]_{t/4T}^{q-1+\delta_{0}} \, dx\\
&& \hspace{-10mm}  \leq (4T)^{\frac{q}{1-\delta_{0}}-1}   \left( \frac{1-\delta_{0}}{q-1+\delta_{0} } \right)   \int_{B_{r_{2}}(z)} [\varphi(|Du|)]^{1-\delta_{0}}  [\varphi(|Du|)]_{t}^{q-1+\delta_{0}} \, dx .
\end{eqnarray}
Likewise, we have
\begin{eqnarray}\label{fubiniinteg3}
\nonumber &&   \hspace{-15mm}  \int_{t_0}^{t} \Lambda^{\frac{q}{1-\delta_{0}}-2}  \int_{\mathcal{E}^{r_{2}}_{\delta_{0}\Lambda/4T} } [ \varphi(|\mathbf{f}|)]^{1-\delta_{0}}   \, dx  \, d\Lambda  \\
&&    \hspace{-5mm}  \leq \left( \frac{4T}{\delta_{0}} \right)^{\frac{q}{1-\delta_{0}}-1}   \left( \frac{1-\delta_{0}}{q-1+\delta_{0} } \right)  \int_{B_{r_{2}}(z)}[\varphi(|\mathbf{f}|)]^{q} \, dx .
\end{eqnarray}
Combining \eqref{fubiniinteg1}, \eqref{fubiniinteg2} and \eqref{fubiniinteg3}, it follows that
 \pagebreak
\begin{eqnarray}\label{fubiniinteg4}
\nonumber &&   \hspace{-10mm} \int_{B_{r_{1}}(z)} [\varphi(|Du|)]^{1-\delta_{0}} [\varphi(|Du|)]_{t}^{q-1+\delta_{0}}\, dx  \\
\nonumber &&  \hspace{-5mm}  \leq T^{\frac{q}{1-\delta_{0}}-1} \Lambda_{0}^{\frac{q}{1-\delta_{0}}}  |B_{r_{2}}(z)|  \\
\nonumber &&   + 320^{n}( \varepsilon T^{\frac{q}{1-\delta_{0}}-1} + c^{*}T^{\frac{q-1-\sigma}{1-\delta_{0}}} )  \int_{B_{r_{2}}(z)} [\varphi(|Du|)]^{1-\delta_{0}}  [\varphi(|Du|)]_{t}^{q-1+\delta_{0}} \, dx\\
\nonumber && \hspace{5mm}  + 320^{n}\frac{ (\varepsilon T^{\frac{q}{1-\delta_{0}}-1} + c^{*}T^{\frac{q-1-\sigma}{1-\delta_{0}}})}{\delta_{0}^{q-1+\delta_{0}} } \int_{B_{r_{2}}(z)}[\varphi(|\mathbf{f}|)]^{q} \, dx  \\
\nonumber &&  \hspace{-5mm}  \leq T^{\frac{q}{1-\delta_{0}}-1} \Lambda_{0}^{\frac{q}{1-\delta_{0}}}  |B_{r_{2}}(z)|  \\
\nonumber &&    + 320^{n}( \varepsilon T^{2} + c^{*}T^{-\frac{\sigma}{2}} )  \int_{B_{r_{2}}(z)} [\varphi(|Du|)]^{1-\delta_{0}}  [\varphi(|Du|)]_{t}^{q-1+\delta_{0}} \, dx\\
&& \hspace{5mm}  + 320^{n}\frac{ ( \varepsilon T^{2} + c^{*}T^{-\frac{\sigma}{2}} )  }{\delta_{0}^{q-1+\delta_{0}} } \int_{B_{r_{2}}(z)}[\varphi(|\mathbf{f}|)]^{q} \, dx,
\end{eqnarray}
since $q \leq 1+\delta_{0} \leq 1+\frac{\sigma}{2}$. Now we choose $T= \left( 2560^n c^{*}\right )^{\frac{2}{\sigma}}$ and $\varepsilon = \frac{1}{2560^n T^{2}}$ in order to have the following estimate.
\begin{eqnarray}\label{fubiniinteg5}
\nonumber &&   \hspace{-10mm} \hspace{0.25em} \Xint-_{B_{r_{1}}(z)} [\varphi(|Du|)]^{1-\delta_{0}} [\varphi(|Du|)]_{t}^{q-1+\delta_{0}}\, dx  \\
\nonumber &&  \hspace{-2.5mm}  \leq \left(  \frac{ c R^n}{(r_{2}-r_{1})^{n}} \hspace{0.25em} \Xint-_{B_{R}(z)}  \left\{ [\varphi(|Du|)]^{1-\delta_{0}}
+ \frac{ [ \varphi(|\mathbf{f}|)]^{1-\delta_{0}}}{\delta_{0}} \right\} \, dx \right)^{\frac{q}{1-\delta_{0}}}   \\
\nonumber &&   \hspace{5mm}  + \frac{1}{2} \hspace{0.25em} \Xint-_{B_{r_{2}}(z)} [\varphi(|Du|)]^{1-\delta_{0}}  [\varphi(|Du|)]_{t}^{q-1+\delta_{0}} \, dx  + c \hspace{0.25em} \Xint-_{B_{R}(z)}[\varphi(|\mathbf{f}|)]^{q} \, dx  .
\end{eqnarray}
Applying Lemma \ref{algeb} with
$$
\phi(s) := \hspace{0.25em} \Xint-_{B_{s}(z)} [\varphi(|Du|)]^{1-\delta_{0}}  [\varphi(|Du|)]_{t}^{q-1+\delta_{0}} \, dx
$$
and $\beta = \frac{qn}{1-\delta_{0}} \in [1,3n]$, we finally obtain
\begin{eqnarray}\label{mainproof1}
\nonumber &&   \hspace{-15mm}  \Xint-_{B_{\frac{R}{2}}(z)} [\varphi(|Du|)]^{1-\delta_{0}} [\varphi(|Du|)]_{t}^{q-1+\delta_{0}}\, dx  \\
 &&  \hspace{-2.5mm}  \leq c \left(   \hspace{0.25em} \Xint-_{B_{R}(z)}   [\varphi(|Du|)]^{1-\delta_{0}}
 \, dx \right)^{\frac{q}{1-\delta_{0}}}   + c \hspace{0.25em} \Xint-_{B_{R}(z)}[\varphi(|\mathbf{f}|)]^{q} \, dx  .
\end{eqnarray}
Letting $t \rightarrow \infty$, we obtain the desired estimate.
\begin{eqnarray}\label{mainproof2}
\notag  &&   \hspace{-15mm} \Xint-_{B_{\frac{R}{2}}(z)}  [\varphi(|Du|)]^{q} \, dx  \\
&& \hspace{-2.5mm}   \leq  c \left[ \left(   \hspace{0.25em} \Xint-_{B_{R}(z)}   [\varphi(|Du|)]^{1-\delta_{0}}
 \, dx \right)^{\frac{q}{1-\delta_{0}}}   + \hspace{0.25em} \Xint-_{B_{R}(z)}[\varphi(|\mathbf{f}|)]^{q} \, dx \right].
\end{eqnarray}
This completes the proof of Theorem \ref{mainthm}.
\end{proof}

\bibliographystyle{amsplain}

\end{document}